# Carleman Linearization of Differential-Algebraic Equations Systems

Marcos A. Hernández-Ortega[a*], C. M. Rergis[a], A. Román-Messina[b], Erlan R. Murillo-Aguirre[b]

[a]*Departamento Académico de Electromecánica, Universidad Autónoma de Guadalajara, Av. Patria 1201, 45129, Zapopan, México*

[b]*Graduate Studies Program in Electrical Engineering, CINVESTAV, Guadalajara, Jal., 45015, Mexico*

[*]Corresponding author. Tel.: +52-33 3648 8824

*E-mail address*: marcos.hernandez@edu.uag.mx (Marcos A. Hernández-Ortega).

**Abstract**

Carleman linearization is a mathematical technique that transforms nonlinear dynamical systems into infinite-dimensional linear systems, enabling simplified analysis. Initially developed for ordinary differential equations (ODEs) and later extended to partial differential equations (PDEs), it has found applications in control theory, biological systems, fluid dynamics, quantum mechanics, finance, and machine learning. This paper extends Carleman linearization to differential-algebraic equation (DAE) systems by introducing auxiliary functions and projecting the resulting system onto a higher-order ODE representation. Theoretical foundations are presented along with conditions under which the transformation is valid. The method is demonstrated on synthetic DAE examples, highlighting its effectiveness even when projection from algebraic variables to state variables is nontrivial or undefined.

*Keywords:* Carleman linearization, differential-algebraic equations systems, nonlinear dynamical systems, ordinary differential equations systems.

## 1. Introduction

In recent decades, the increasing prevalence of nonlinear dynamical systems in science and technology has underscored the need for robust analytical methodologies. Complex systems such as biological networks, smart energy grids, and advanced transportation infrastructures exhibit nonlinear behavior, high interdependence, and evolving dynamics that often defy traditional modeling approaches [1]-[3]. As these systems grow in sophistication, effective tools are needed for prediction, control design, optimization, and monitoring [4]-[5]. A key strategy in managing this complexity is modular modeling, which decomposes systems into smaller, reusable subsystems. Representing such models as graphs—where modules are nodes and interconnections are edges—enhances clarity and scalability, especially in fields such as biology, engineering, and computer science [6]-[8].

To simulate and analyze these modular, interconnected models, Differential-Algebraic Equation (DAE) frameworks have become foundational. By combining differential and algebraic constraints, DAEs effectively capture both temporal evolution and structural properties. Their applications span vehicle dynamics, robotics, electrical circuits, and multiphysics simulations, and their use is deeply integrated into modeling languages. DAEs also play a growing role in machine learning and neural-based simulation frameworks for large-scale, constrained systems [6]-[8]. DAE models are valued for their ability to accurately represent the behavior of complex systems in engineering, physics, and cyber-physical domains, and they are gaining influence in advanced machine learning methodologies as well [9]. This cements their role as an essential tool for the study of modern dynamical systems.

One of the key strengths of DAE-based modelling lies in its ability to support modularity and component reuse, which is particularly beneficial in multiphysics simulations. This feature is integral to modelling environments such as Modelica, where DAE formulations naturally arise due to component interconnections and algebraic constraints [9].

Differential-Algebraic Equation (DAE) systems are widely used in scientific and engineering domains, particularly for modeling vehicle dynamics, trajectory planning in mobile robotics, and advanced control systems [10]. Their versatility extends to the analysis of mechanical multibody systems, electrical circuits, and fluid dynamics problems, where they provide a unified framework for capturing both dynamic behavior and algebraic constraints [11].

In the power systems domain, DAEs are crucial for capturing the interactions between electrical network components and dynamic elements, such as generators and loads [12]. More broadly, DAE frameworks range from classical formulations of constrained mechanical systems to modern implementations involving machine learning and neural networks, which can efficiently simulate large-scale dynamic systems under multiple constraints [13].

Recent work has demonstrated the efficacy of DAE-based models in electric vehicle simulation, further emphasizing their growing relevance in cyber-physical and energy systems research [14]. In parallel, advances in computer science and numerical methods continue to drive the development of more robust and scalable algorithms for solving DAEs, reinforcing their foundational role in modern system modelling and control [15].

One methodological challenge lies in analyzing nonlinear systems represented by DAEs. Carleman linearization, a classical technique that converts finite-dimensional nonlinear systems into infinite-dimensional linear ones, offers a promising solution. It does so by expanding nonlinear terms into higher-order monomials, which are introduced as new variables, yielding an infinite set of linear equations. The core idea involves expressing nonlinear terms as higher-order monomials and introducing them as additional variables.

The result is a linear system of infinite dimension that approximates the dynamics of the original nonlinear model [16].

Carleman linearization has demonstrated success across diverse domains, including nonlinear circuit analysis, chemical kinetics for quantum computing, and power systems [16-23]. In circuit theory, it has been employed to analyze and design nonlinear electronic systems [17]. In chemical kinetics, it enables the transformation of nonlinear ordinary differential equations (ODEs) into linear forms, which facilitates simulation on quantum computing platforms [18]. The technique has also been used to develop higher-order small-signal models in power system analysis [19].

Moreover, Carleman linearization supports the systematic expansion of certain partial differential equation (PDE) systems [20]. It has also proven effective for reachability analysis in weakly nonlinear dynamical systems, offering new avenues for verifying system behavior [21].

Of particular interest, Carleman linearization has facilitated the application of advanced nonlinear analysis techniques such as the method of normal forms [22] and the perturbed Koopman mode analysis [23]. While efforts have been made to extend similar methodologies to DAE systems using singular perturbation analysis [12], such approaches are sensitive to the choice of perturbation parameters and may introduce spurious dynamics, thereby limiting their reliability.

This research introduces a novel extension of Carleman linearization tailored specifically for nonlinear DAE systems. By augmenting the system with additional functions that capture nonlinear interactions, the method projects the DAE into an equivalent higher-order linear ODE model. This transformation supports systematic analysis and control design, even when explicit algebraic-to-state mappings are unavailable.

The theoretical foundations of the method are thoroughly developed, and its practical usefulness is demonstrated through synthetic DAE case studies. Results confirm the method's effectiveness in approximating DAE dynamics, with implications for model order reduction, stability analysis, and advanced control design. This work opens the door for using nonlinear analysis techniques—such as normal forms and Koopman-based methods—in DAE frameworks.

## 2. Carleman Linearization for ODE Systems

Consider a nonlinear dynamic system described by the ordinary differential equation:

$$\dot{\mathbf{x}} = f(\mathbf{x}, \mathbf{u}) \tag{1}$$

where $\mathbf{x} \in \mathcal{R}^N$ is the vector of state variables, $\mathbf{u} \in \mathcal{R}^P$ is the vector of system inputs, and $f$ is a function describing the system evolution in time and is holomorphic around a stable equilibrium point (*sep*) given by the vectors $\mathbf{x}_{sep}$, $\mathbf{u}_{sep}$.

Assume that $f$ is holomorphic in a neighborhood of a *sep*. Under small perturbations in the state or input vectors, the nonlinear system (1) can be approximated using a Taylor series expansion around this equilibrium point as:

$$\dot{\mathbf{x}} \approx \Delta\dot{\mathbf{x}} = F_1(\Delta\mathbf{x}, \Delta\mathbf{u}) + F_2(\Delta\mathbf{x}, \Delta\mathbf{u}) + F_3(\Delta\mathbf{x}, \Delta\mathbf{u}) + h.o.t. \tag{2}$$

where $F_j(\Delta\mathbf{x}, \Delta\mathbf{u})$ denotes the $j$-th term of the Taylor series and higher order terms (*h.o.t.*) accounts for all terms beyond a chosen truncation level.

For simplicity, this article considers only perturbations of the state variables, assuming the inputs remain constant. Consequently, the expression in (2) simplifies to the following truncated form [22]:

$$\Delta\dot{\mathbf{x}} = \mathbf{A}_{1,1}\Delta\mathbf{x} + \mathbf{A}_{1,2}\Delta\mathbf{x}^{[2]} + \mathbf{A}_{1,3}\Delta\mathbf{x}^{[3]} + h.o.t. \tag{3}$$

where the Kronecker power of the perturbation vector is defined as:

$$\Delta\mathbf{x}^{[n]} = \underbrace{\Delta\mathbf{x} \otimes \Delta\mathbf{x} \otimes \cdots \otimes \Delta\mathbf{x}}_{n\ terms}. \tag{4}$$

In expression (3), the symbol $\otimes$ represents the Kronecker product. The matrix $\mathbf{A}_{1,1}$ represents the Jacobian of the system evaluated at the equilibrium point, while $\mathbf{A}_{1,2}$ and $\mathbf{A}_{1,3}$ correspond to the second and third order coefficient matrices, respectively [24].

Subsequently, the higher-order variables $\Delta \mathbf{x}^{[n]}$ are introduced as extended states, and their time derivatives are defined as follows:

$$\Delta \dot{\mathbf{x}}^{[n]} = \underbrace{\Delta \dot{\mathbf{x}} \otimes \Delta \mathbf{x} \otimes \cdots \otimes \Delta \mathbf{x}}_{n\ terms} + \underbrace{\Delta \mathbf{x} \otimes \Delta \dot{\mathbf{x}} \otimes \cdots \otimes \Delta \mathbf{x}}_{n\ terms} + \cdots + \underbrace{\Delta \mathbf{x} \otimes \Delta \mathbf{x} \otimes \cdots \otimes \Delta \dot{\mathbf{x}}}_{n\ terms}. \quad (5)$$

Incorporating these terms into the state-space representation of the system yields the following block-structured dynamic model:

$$\underbrace{\begin{bmatrix} \Delta \dot{\mathbf{x}} \\ \Delta \dot{\mathbf{x}}^{[2]} \\ \Delta \dot{\mathbf{x}}^{[3]} \\ \vdots \end{bmatrix}}_{\Delta \dot{\mathbf{x}}_{nord}} = \underbrace{\begin{bmatrix} \mathbf{A}_{1,1} & \mathbf{A}_{1,2} & \mathbf{A}_{1,3} & \cdots \\ 0 & \mathbf{A}_{2,2} & \mathbf{A}_{2,3} & \cdots \\ 0 & 0 & \mathbf{A}_{3,3} & \cdots \\ \vdots & \vdots & \vdots & \ddots \end{bmatrix}}_{\mathbf{A}_{nord}} \underbrace{\begin{bmatrix} \Delta \mathbf{x} \\ \Delta \mathbf{x}^{[2]} \\ \Delta \mathbf{x}^{[3]} \\ \vdots \end{bmatrix}}_{\Delta \mathbf{x}_{nord}} \quad (6)$$

where the block matrices $\mathbf{A}_{i,j}$ are recursively defined as:

$$\mathbf{A}_{i,j} = \mathbf{A}_{1,j} \otimes \mathbf{I}_{n^{i-1}} + \mathbf{I}_n \otimes \mathbf{A}_{i-1,j} \quad (7)$$

for $i = 2, 3, \ldots, l$ and $j = 1, 2, \ldots, l - i + 1$ and where $\mathbf{I}_n$ denotes the $n \times n$ identity matrix [22]. Thus, for $i > 2$, the matrices $\mathbf{A}_{i,j}$ are recursively constructed based on the initial matrices $\mathbf{A}_{1,j}$.

As previously discussed, the Carleman linearization technique enables the transformation of a nonlinear ODE system of the form Eq. (1) into an infinite-dimensional linear system, as represented in Eq. (6). This linear system can be truncated to a finite order to achieve an approximate representation of the original nonlinear dynamics. A notable property of the Carleman-linearized system (Eq. 6) is that its eigenvalues are composed of linear combinations of the eigenvalues of the Jacobian matrix $\mathbf{A}_{1,1}$, denoted as $\lambda_i$. Specifically, the spectrum of the linearized system (6) includes terms such as $\lambda_i$, $\lambda_i + \lambda_j$, $\lambda_i + \lambda_j + \lambda_k$, and so on.

In the following section, this linearization methodology is extended to differential-algebraic equation (DAE) systems.

## 3. Carleman Linearization of DAE Systems

*3.1 Background*

Consider a nonlinear DAE system of the form:

$$\dot{\mathbf{x}} = \mathbf{g}(\mathbf{x}, \mathbf{z}, \mathbf{u}) \tag{8}$$

$$\mathbf{0} = \mathbf{h}(\mathbf{x}, \mathbf{z}, \mathbf{u}) \tag{9}$$

where $\mathbf{z} \in \mathcal{R}^M$ is a vector of algebraic variables. Around the operating point ($\mathbf{x}_{sep}$, $\mathbf{z}_{sep}$, $\mathbf{u}_{sep}$), there exists a unique value of $\mathbf{z}$ that satisfies equation (9) for each pair ($\mathbf{x}, \mathbf{u}$).

In equations (8), (9) the terms $\mathbf{g}, \mathbf{h}$ are nonlinear and describe the time evolution and algebraic constraints of the system. Both functions are assumed to be holomorphic in a neighborhood of a stable equilibrium point defined by the vectors $\mathbf{x}_{sep}$, $\mathbf{z}_{sep}$, $\mathbf{u}_{sep}$. Moreover, the Jacobian matrix of $\mathbf{h}$ with respect to the variables evaluated at the *sep* and denoted $\mathbf{H}_{1,4}$, is invertible; that is:

$$|\mathbf{H}_{1,4}| = \left| \begin{bmatrix} \frac{\partial h_1}{\partial z_1}\Big|_{\mathbf{x}_{sep},\mathbf{z}_{sep},\mathbf{u}_{sep}} & \cdots & \frac{\partial h_1}{\partial z_M}\Big|_{\mathbf{x}_{sep},\mathbf{z}_{sep},\mathbf{u}_{sep}} \\ \vdots & \ddots & \vdots \\ \frac{\partial h_M}{\partial z_1}\Big|_{\mathbf{x}_{sep},\mathbf{z}_{sep},\mathbf{u}_{sep}} & \cdots & \frac{\partial h_M}{\partial z_M}\Big|_{\mathbf{x}_{sep},\mathbf{z}_{sep},\mathbf{u}_{sep}} \end{bmatrix} \right| \neq 0 \tag{10}$$

Additionally, for every value of $\mathbf{x}$ in the neighborhood of $\mathbf{x}_{sep}$, there is a unique value of $\mathbf{z}$ near $\mathbf{z}_{sep}$ that satisfies (9).

*3.2 Taylor Series Approximation*

In order to analyze the nonlinear dynamics of the DAE system defined by Eqs. (8) and (9), a Taylor series expansion can be employed, following an approach analogous to that used for the ODE system in Eq. (2), leading to:

$$\begin{bmatrix} \dot{\mathbf{x}} \\ \mathbf{0} \end{bmatrix} \approx \begin{bmatrix} \Delta\dot{\mathbf{x}} \\ \mathbf{0} \end{bmatrix} = \begin{bmatrix} G_1(\Delta\mathbf{x}, \Delta\mathbf{u}) \\ H_1(\Delta\mathbf{x}, \Delta\mathbf{u}) \end{bmatrix} + \begin{bmatrix} G_2(\Delta\mathbf{x}, \Delta\mathbf{u}) \\ H_2(\Delta\mathbf{x}, \Delta\mathbf{u}) \end{bmatrix} + \begin{bmatrix} G_3(\Delta\mathbf{x}, \Delta\mathbf{u}) \\ H_3(\Delta\mathbf{x}, \Delta\mathbf{u}) \end{bmatrix} + h.o.t. \tag{11}$$

where $G_j(\Delta\mathbf{x}, \Delta\mathbf{u})$ and $H_j(\Delta\mathbf{x}, \Delta\mathbf{u})$ denote the *j*-th terms of the Taylor series expansion related to the functions $g$ and $h$, respectively.

By neglecting the input variables for simplicity, equations (12) and (13) can be derived through the expansion of Eq. (11).

$$\Delta\dot{\mathbf{x}} = \mathbf{G}_{1,1}\Delta\mathbf{x} + \mathbf{G}_{1,2}\Delta\mathbf{x}^{[2]} + \mathbf{G}_{1,3}\Delta\mathbf{x}^{[3]} + \mathbf{G}_{1,4}\Delta\mathbf{z} + \mathbf{G}_{1,5}(\Delta\mathbf{x}\otimes\Delta\mathbf{z}) + \\ +\mathbf{G}_{1,6}\Delta\mathbf{z}^{[2]} + \mathbf{G}_{1,7}(\Delta\mathbf{x}^{[2]}\otimes\Delta\mathbf{z}) + \mathbf{G}_{1,8}(\Delta\mathbf{x}\otimes\Delta\mathbf{z}^{[2]}) + \mathbf{G}_{1,9}\Delta\mathbf{z}^{[3]} + h.o.t. \quad (12)$$

$$\mathbf{0} = \Delta\mathbf{h} = \mathbf{H}_{1,1}\Delta\mathbf{x} + \mathbf{H}_{1,2}\Delta\mathbf{x}^{[2]} + \mathbf{H}_{1,3}\Delta\mathbf{x}^{[3]} + \mathbf{H}_{1,4}\Delta\mathbf{z} + \mathbf{H}_{1,5}(\Delta\mathbf{x}\otimes\Delta\mathbf{z}) + \\ +\mathbf{H}_{1,6}\Delta\mathbf{z}^{[2]} + \mathbf{H}_{1,7}(\Delta\mathbf{x}^{[2]}\otimes\Delta\mathbf{z}) + \mathbf{H}_{1,8}(\Delta\mathbf{x}\otimes\Delta\mathbf{z}^{[2]}) + \mathbf{H}_{1,9}\Delta\mathbf{z}^{[3]} + h.o.t. \quad (13)$$

In expressions (12) and (13), the terms $\Delta\mathbf{z}^{[k]}$ are defined analogously to the definition given in (4). The analysis primarily focuses on terms of third order and lower. According to the Carleman linearization approach, the time derivatives $\Delta\dot{\mathbf{x}}^{[2]}$ can be obtained using equation (5). Subsequently, by substituting Eq. (12) and applying the properties of the Kronecker product [24], the following expressions are derived:

$$\Delta\dot{\mathbf{x}}^{[2]} = \Delta\dot{\mathbf{x}}\otimes\Delta\mathbf{x} + \Delta\mathbf{x}\otimes\Delta\dot{\mathbf{x}}$$

$$\Delta\dot{\mathbf{x}}^{[2]} = \big(\mathbf{G}_{1,1}\Delta\mathbf{x} + \mathbf{G}_{1,2}\Delta\mathbf{x}^{[2]} + \mathbf{G}_{1,3}\Delta\mathbf{x}^{[3]} + \mathbf{G}_{1,4}\Delta\mathbf{z} + \mathbf{G}_{1,5}(\Delta\mathbf{x}\otimes\Delta\mathbf{z}) + \\ +\mathbf{G}_{1,6}\Delta\mathbf{z}^{[2]} + \mathbf{G}_{1,7}(\Delta\mathbf{x}^{[2]}\otimes\Delta\mathbf{z}) + \mathbf{G}_{1,8}(\Delta\mathbf{x}\otimes\Delta\mathbf{z}^{[2]}) + \mathbf{G}_{1,9}\Delta\mathbf{z}^{[3]} + h.o.t.\big)\otimes\Delta\mathbf{x} + \\ \Delta\mathbf{x}\otimes\big(\mathbf{G}_{1,1}\Delta\mathbf{x} + \mathbf{G}_{1,2}\Delta\mathbf{x}^{[2]} + \mathbf{G}_{1,3}\Delta\mathbf{x}^{[3]} + \mathbf{G}_{1,4}\Delta\mathbf{z} + \mathbf{G}_{1,5}(\Delta\mathbf{x}\otimes\Delta\mathbf{z}) + \\ +\mathbf{G}_{1,6}\Delta\mathbf{z}^{[2]} + \mathbf{G}_{1,7}(\Delta\mathbf{x}^{[2]}\otimes\Delta\mathbf{z}) + \mathbf{G}_{1,8}(\Delta\mathbf{x}\otimes\Delta\mathbf{z}^{[2]}) + \mathbf{G}_{1,9}\Delta\mathbf{z}^{[3]} + h.o.t.\big)$$

$$\Delta\dot{\mathbf{x}}^{[2]} = (\mathbf{G}_{1,1}\otimes\mathbf{I}_N)\Delta\mathbf{x}^{[2]} + (\mathbf{G}_{1,2}\otimes\mathbf{I}_{N^2})\Delta\mathbf{x}^{[3]} + (\mathbf{G}_{1,3}\otimes\mathbf{I}_{N^3})\Delta\mathbf{x}^{[4]} + \\ +(\mathbf{G}_{1,4}\otimes\mathbf{I}_N)(\Delta\mathbf{z}\otimes\Delta\mathbf{x}) + (\mathbf{G}_{1,5}\otimes\mathbf{I}_N)(\Delta\mathbf{x}\otimes\Delta\mathbf{z}\otimes\Delta\mathbf{x}) + (\mathbf{G}_{1,6}\otimes\mathbf{I}_N)(\Delta\mathbf{z}^{[2]}\otimes\Delta\mathbf{x}) + \\ +(\mathbf{G}_{1,7}\otimes\mathbf{I}_N)(\Delta\mathbf{x}^{[2]}\otimes\Delta\mathbf{z}\otimes\Delta\mathbf{x}) + (\mathbf{G}_{1,8}\otimes\mathbf{I}_N)(\Delta\mathbf{x}\otimes\Delta\mathbf{z}^{[2]}\otimes\Delta\mathbf{x}) + (\mathbf{G}_{1,9}\otimes\mathbf{I}_N)(\Delta\mathbf{z}^{[3]}\otimes\Delta\mathbf{x}) + \\ +(\mathbf{I}_N\otimes\mathbf{G}_{1,1})\Delta\mathbf{x}^{[2]} + (\mathbf{I}_N\otimes\mathbf{G}_{1,2})\Delta\mathbf{x}^{[3]} + (\mathbf{I}_N\otimes\mathbf{G}_{1,3})\Delta\mathbf{x}^{[4]} + (\mathbf{I}_N\otimes\mathbf{G}_{1,4})(\Delta\mathbf{x}\otimes\Delta\mathbf{z}) + \\ +(\mathbf{I}_N\otimes\mathbf{G}_{1,5})(\Delta\mathbf{x}^{[2]}\otimes\Delta\mathbf{z}) + (\mathbf{I}_N\otimes\mathbf{G}_{1,6})(\Delta\mathbf{x}\otimes\Delta\mathbf{z}^{[2]}) + (\mathbf{I}_N\otimes\mathbf{G}_{1,7})(\Delta\mathbf{x}^{[3]}\otimes\Delta\mathbf{z}) + \\ +(\mathbf{I}_N\otimes\mathbf{G}_{1,8})(\Delta\mathbf{x}^{[2]}\otimes\Delta\mathbf{z}^{[2]}) + (\mathbf{I}_N\otimes\mathbf{G}_{1,9})(\Delta\mathbf{x}\otimes\Delta\mathbf{z}^{[3]}) + h.o.t.$$

(14)

For simplicity, terms of order four and higher are omitted. By defining the matrices $\mathbf{G}_{2,2}$ and $\mathbf{G}_{2,3}$ analogously to definition (7), the expression (14) can be simplified as follows:

$$\Delta\dot{\mathbf{x}}^{[2]} = \mathbf{G}_{2,2}\Delta\mathbf{x}^{[2]} + \mathbf{G}_{2,3}\Delta\mathbf{x}^{[3]} + (\mathbf{G}_{1,4}\otimes\mathbf{I}_N)(\Delta\mathbf{z}\otimes\Delta\mathbf{x}) +$$
$$+ (\mathbf{G}_{1,5}\otimes\mathbf{I}_N)(\Delta\mathbf{x}\otimes\Delta\mathbf{z}\otimes\Delta\mathbf{x}) + (\mathbf{G}_{1,6}\otimes\mathbf{I}_N)(\Delta\mathbf{z}^{[2]}\otimes\Delta\mathbf{x}) +$$
$$+ (\mathbf{I}_N\otimes\mathbf{G}_{1,4})(\Delta\mathbf{x}\otimes\Delta\mathbf{z}) + (\mathbf{I}_N\otimes\mathbf{G}_{1,5})(\Delta\mathbf{x}^{[2]}\otimes\Delta\mathbf{z}) +$$
$$+ (\mathbf{I}_N\otimes\mathbf{G}_{1,6})(\Delta\mathbf{x}\otimes\Delta\mathbf{z}^{[2]}) + h.o.t. \quad (15)$$

Equation (15) contains structurally equivalent terms such as $\Delta\mathbf{z}\otimes\Delta\mathbf{x}$ and $\Delta\mathbf{x}\otimes\Delta\mathbf{z}$. To consolidate the representation in equation (14), a set of matrices $\widetilde{\mathbf{G}}_{1,j}$ is introduced. These matrices contain the same information as $(\mathbf{G}_{1,j}\otimes\mathbf{I}_N)$, but are reordered to align with the current term structure. The result of this reformulation is presented in equation (16).

$$\Delta\dot{\mathbf{x}}^{[2]} = \mathbf{G}_{2,2}\Delta\mathbf{x}^{[2]} + \mathbf{G}_{2,3}\Delta\mathbf{x}^{[3]} + (\widetilde{\mathbf{G}}_{1,4} + \mathbf{I}_N\otimes\mathbf{G}_{1,4})(\Delta\mathbf{x}\otimes\Delta\mathbf{z}) +$$
$$+ (\widetilde{\mathbf{G}}_{1,5} + \mathbf{I}_N\otimes\mathbf{G}_{1,5})(\Delta\mathbf{x}^{[2]}\otimes\Delta\mathbf{z}) + (\widetilde{\mathbf{G}}_{1,6} + \mathbf{I}_N\otimes\mathbf{G}_{1,6})(\Delta\mathbf{x}\otimes\Delta\mathbf{z}^{[2]}) + h.o.t. \quad (16)$$

Then, by defining $\mathbf{G}_{2,5} = \widetilde{\mathbf{G}}_{1,4} + \mathbf{I}_N\otimes\mathbf{G}_{1,4}$, $\mathbf{G}_{2,7} = \widetilde{\mathbf{G}}_{1,5} + \mathbf{I}_N\otimes\mathbf{G}_{1,5}$, and $\mathbf{G}_{2,8} = \widetilde{\mathbf{G}}_{1,6} + \mathbf{I}_N\otimes\mathbf{G}_{1,6}$, the expression (16) can be rewritten in the following form:

$$\Delta\dot{\mathbf{x}}^{[2]} = \mathbf{G}_{2,2}\Delta\mathbf{x}^{[2]} + \mathbf{G}_{2,3}\Delta\mathbf{x}^{[3]} + \mathbf{G}_{2,5}(\Delta\mathbf{x}\otimes\Delta\mathbf{z}) +$$
$$+ \mathbf{G}_{2,7}(\Delta\mathbf{x}^{[2]}\otimes\Delta\mathbf{z}) + \mathbf{G}_{2,8}(\Delta\mathbf{x}\otimes\Delta\mathbf{z}^{[2]}) + h.o.t. \quad (17)$$

Now, let us develop the time derivative $\Delta\dot{\mathbf{x}}^{[3]}$. In this case, the fourth- and higher-order terms are ignored since the first steps for simplicity.

$$\Delta\dot{\mathbf{x}}^{[3]} = \Delta\dot{\mathbf{x}}\otimes\Delta\mathbf{x}^{[2]} + \Delta\mathbf{x}\otimes\Delta\dot{\mathbf{x}}\otimes\Delta\mathbf{x} + \Delta\mathbf{x}^{[2]}\otimes\Delta\dot{\mathbf{x}}$$

$$\Delta\dot{\mathbf{x}}^{[3]} = (\mathbf{G}_{1,1}\Delta\mathbf{x} + \mathbf{G}_{1,2}\Delta\mathbf{x}^{[2]} + \mathbf{G}_{1,3}\Delta\mathbf{x}^{[3]} + \mathbf{G}_{1,4}\Delta\mathbf{z} + \mathbf{G}_{1,5}(\Delta\mathbf{x}\otimes\Delta\mathbf{z}) +$$
$$+ \mathbf{G}_{1,6}\Delta\mathbf{z}^{[2]} + \mathbf{G}_{1,7}(\Delta\mathbf{x}^{[2]}\otimes\Delta\mathbf{z}) + \mathbf{G}_{1,8}(\Delta\mathbf{x}\otimes\Delta\mathbf{z}^{[2]}) + \mathbf{G}_{1,9}\Delta\mathbf{z}^{[3]} + h.o.t.)\otimes\Delta\mathbf{x}^{[2]} +$$
$$+ \Delta\mathbf{x}\otimes(\mathbf{G}_{1,1}\Delta\mathbf{x} + \mathbf{G}_{1,2}\Delta\mathbf{x}^{[2]} + \mathbf{G}_{1,3}\Delta\mathbf{x}^{[3]} + \mathbf{G}_{1,4}\Delta\mathbf{z} + \mathbf{G}_{1,5}(\Delta\mathbf{x}\otimes\Delta\mathbf{z}) +$$
$$+ \mathbf{G}_{1,6}\Delta\mathbf{z}^{[2]} + \mathbf{G}_{1,7}(\Delta\mathbf{x}^{[2]}\otimes\Delta\mathbf{z}) + \mathbf{G}_{1,8}(\Delta\mathbf{x}\otimes\Delta\mathbf{z}^{[2]}) + \mathbf{G}_{1,9}\Delta\mathbf{z}^{[3]} + h.o.t.)\otimes\Delta\mathbf{x} +$$
$$\Delta\mathbf{x}^{[2]}\otimes(\mathbf{G}_{1,1}\Delta\mathbf{x} + \mathbf{G}_{1,2}\Delta\mathbf{x}^{[2]} + \mathbf{G}_{1,3}\Delta\mathbf{x}^{[3]} + \mathbf{G}_{1,4}\Delta\mathbf{z} + \mathbf{G}_{1,5}(\Delta\mathbf{x}\otimes\Delta\mathbf{z}) +$$
$$+ \mathbf{G}_{1,6}\Delta\mathbf{z}^{[2]} + \mathbf{G}_{1,7}(\Delta\mathbf{x}^{[2]}\otimes\Delta\mathbf{z}) + \mathbf{G}_{1,8}(\Delta\mathbf{x}\otimes\Delta\mathbf{z}^{[2]}) + \mathbf{G}_{1,9}\Delta\mathbf{z}^{[3]} + h.o.t.)$$

$$\Delta\dot{\mathbf{x}}^{[3]} = (\mathbf{G}_{1,1}\otimes\mathbf{I}_{N^2})\Delta\mathbf{x}^{[3]} + (\mathbf{G}_{1,4}\otimes\mathbf{I}_{N^2})(\Delta\mathbf{z}\otimes\Delta\mathbf{x}^{[2]}) +$$
$$+ (\mathbf{I}_N\otimes\mathbf{G}_{1,1}\otimes\mathbf{I}_N)\Delta\mathbf{x}^{[3]} + (\mathbf{I}_N\otimes\mathbf{G}_{1,4}\otimes\mathbf{I}_N)(\Delta\mathbf{x}\otimes\Delta\mathbf{z}\otimes\Delta\mathbf{x}) + \quad (18)$$
$$+ (\mathbf{I}_{N^2}\otimes\mathbf{G}_{1,1})\Delta\mathbf{x}^{[3]} + (\mathbf{I}_{N^2}\otimes\mathbf{G}_{1,4})(\Delta\mathbf{x}^{[2]}\otimes\Delta\mathbf{z}) + h.o.t.$$

Again, we define matrix $\mathbf{G}_{3,3}$ according to (7) as $\mathbf{G}_{3,3} = \mathbf{G}_{1,1}\otimes\mathbf{I}_{N^2} + \mathbf{I}_N\otimes\mathbf{G}_{1,1}\otimes\mathbf{I}_N + \mathbf{I}_{N^2}\otimes\mathbf{G}_{1,1}$. Additionally, by introducing the matrices $\widetilde{\mathbf{G}}_{1,4}$ and $\widehat{\mathbf{G}}_{1,4}$ which contain the same information as the terms $\mathbf{G}_{1,4} \otimes \mathbf{I}_{N^2}$ and $\mathbf{I}_N \otimes \mathbf{G}_{1,4} \otimes \mathbf{I}_N$, respectively, but reordered to match the current term structure. Using these definitions, $\mathbf{G}_{3,7}$ can be defined as $\widetilde{\mathbf{G}}_{1,4} + \widehat{\mathbf{G}}_{1,4} + \mathbf{I}_{N^2} \otimes \mathbf{G}_{1,4}$. Consequently, equation (18) can be rewritten as follows:

$$\Delta\dot{\mathbf{x}}^{[3]} = \mathbf{G}_{3,3}\Delta\mathbf{x}^{[3]} + \mathbf{G}_{3,7}(\Delta\mathbf{x}^{[2]} \otimes \Delta\mathbf{z}) + h.o.t. \tag{19}$$

Using expressions (12), (13), (17), and (19), the Carleman extended model of the DAE system is obtained as follows:

$$\begin{bmatrix} \Delta\dot{\mathbf{x}} \\ \Delta\dot{\mathbf{x}}^{[2]} \\ \Delta\dot{\mathbf{x}}^{[3]} \\ 0 \\ \vdots \end{bmatrix} = \begin{bmatrix} \mathbf{G}_{1,1} & \mathbf{G}_{1,2} & \mathbf{G}_{1,3} & \mathbf{G}_{1,4} & \mathbf{G}_{1,5} & \mathbf{G}_{1,6} & \mathbf{G}_{1,7} & \mathbf{G}_{1,8} & \mathbf{G}_{1,9} \\ 0 & \mathbf{G}_{2,2} & \mathbf{G}_{2,3} & 0 & \mathbf{G}_{2,5} & 0 & \mathbf{G}_{2,7} & \mathbf{G}_{2,8} & 0 \\ 0 & 0 & \mathbf{G}_{3,3} & 0 & 0 & 0 & \mathbf{G}_{3,7} & 0 & 0 \\ \mathbf{H}_{1,1} & \mathbf{H}_{1,2} & \mathbf{H}_{1,3} & \mathbf{H}_{1,4} & \mathbf{H}_{1,5} & \mathbf{H}_{1,6} & \mathbf{H}_{1,7} & \mathbf{H}_{1,8} & \mathbf{H}_{1,9} \\ \vdots & \vdots & \vdots & \vdots & \vdots & \vdots & \vdots & \vdots & \vdots \end{bmatrix} \begin{bmatrix} \Delta\mathbf{x} \\ \Delta\mathbf{x}^{[2]} \\ \Delta\mathbf{x}^{[3]} \\ \Delta\mathbf{z} \\ \Delta\mathbf{x} \otimes \Delta\mathbf{z} \\ \Delta\mathbf{z}^{[2]} \\ \Delta\mathbf{x}^{[2]} \otimes \Delta\mathbf{z} \\ \Delta\mathbf{x} \otimes \Delta\mathbf{z}^{[2]} \\ \Delta\mathbf{z}^{[3]} \end{bmatrix} \tag{20}$$

$$+h.o.t.$$

As it can be noted, if the Carleman linearization is made up to a certain order of approximation equal or greater than two, the resulting system is non-determined. This makes the transformation of the system (20) into an equivalent ODE representation difficult, to which several nonlinear analysis techniques, such as the method of normal forms [22] and the perturbed Koopman mode analysis [23], can be applied. If only the first-order terms are considered, a Kron reduction can be easily applied [12].

Therefore, a new method is sought that can transform the system (20) into an equivalent higher-order linear ODE representation. Some useful developments are presented next.

*3.3 Equivalent ODE Representation with an Implicit Function*

<u>Proposition</u> 1. Consider a nonlinear DAE system of the form (8), (9) with the conditions described in Section 3.1. Then, there exists a function $\widetilde{\mathbf{h}}$ such that

$$\mathbf{z} = \widetilde{\mathbf{h}}(\mathbf{x}, \mathbf{u}) \tag{21}$$

*Proof*:

The existence of the function $\tilde{h}$ follows from the Theorem of the Implicit Function and the uniqueness of $\mathbf{z}$ for each value of $\mathbf{x}$ about a *sep*. According to the Implicit Function Theorem, a function $\tilde{h}$ exists if: (i) The function $h$ in (9) is differentiable about a certain point $(\mathbf{x}_0, \mathbf{z}_0, \mathbf{u}_0)$, and (ii) The Jacobian matrix of the function $h$ about $(\mathbf{x}_0, \mathbf{z}_0, \mathbf{u}_0)$ is non-singular (i.e., its determinant is non-zero) [25].

If we consider the point $(\mathbf{x}_0, \mathbf{z}_0, \mathbf{u}_0)$ to be a *sep* defined by $(\mathbf{x}_{sep}, \mathbf{z}_{sep}, \mathbf{u}_{sep})$, as defined earlier, then both conditions of the Implicit Function Theorem are satisfied. Furthermore, the uniqueness of $\mathbf{z}$ with respect to $\mathbf{x}$ ensures the existence of $\tilde{h}$ for the DAE system as described in Section 3.1.

∎

Theorem 1. Consider a nonlinear DAE system of the form (8), (9). If there exists a function $\tilde{h}$ of the form (21), then, a higher-order linear ODE representation of the form (6) can be derived.

*Proof*:

If the function $\tilde{h}$ defined in (21) exists, then a Taylor series approximation can be derived, yielding the following expression:

$$\Delta \mathbf{z} = \widetilde{\mathbf{H}}_{1,1} \Delta \mathbf{x} + \widetilde{\mathbf{H}}_{1,2} \Delta \mathbf{x}^{[2]} + \widetilde{\mathbf{H}}_{1,3} \Delta \mathbf{x}^{[3]} + h.o.t. \qquad (22)$$

Substituting expression (22) into equation (20) and performing standard algebraic manipulations yields an equivalent system of the form (6). However, although $\tilde{h}$ may exist, obtaining it for many dynamical systems is often a difficult and impractical task.

∎

Consequently, if the function $\tilde{h}$ is unavailable, the challenge shifts to determining the matrices $\widetilde{\mathbf{H}}_{1,j}$ on the right-hand part of equation (22).

*3.4 Procedure for Obtaining Matrices $\widetilde{\mathbf{H}}_{j,1}$*

To determine the matrices $\widetilde{\mathbf{H}}_{1,j}$, we begin by revisiting the undetermined higher-order system presented in equation (13). To render system (20) square, additional terms are introduced into the formulation. In the Carleman linearization method for ODE systems, the terms $\Delta \mathbf{x}^{[k]}$ are introduced as additional state variables. Their corresponding time derivatives are then computed and incorporated into an extended linear system. However, unlike in ODE systems, the terms $\Delta \mathbf{z}^{[k]}$ and $\Delta \mathbf{x}^{[j]} \otimes \Delta \mathbf{z}^{[k]}$ in DAE systems cannot be treated as state variables. Instead, the use of the terms $\mathbf{0} = \Delta \boldsymbol{h}^{[j]}$ and $\mathbf{0} = \Delta \boldsymbol{h}^{[j]} \otimes \Delta \mathbf{x}^{[k]}$ is proposed, where $\Delta \boldsymbol{h}^{[j]}$ is defined as follows:

$$\mathbf{0} = \Delta \boldsymbol{h}^{[j]} = \underbrace{\Delta \boldsymbol{h} \otimes \Delta \boldsymbol{h} \otimes \cdots \otimes \Delta \boldsymbol{h}}_{j \; terms}. \tag{23}$$

Expanding the second-order term $\Delta \boldsymbol{h}^{[2]}$ yields:

$$\Delta \boldsymbol{h}^{[2]} = \big(\mathbf{H}_{1,1}\Delta\mathbf{x} + \mathbf{H}_{1,2}\Delta\mathbf{x}^{[2]} + \mathbf{H}_{1,4}\Delta\mathbf{z} + \mathbf{H}_{1,5}(\Delta\mathbf{x}\otimes\Delta\mathbf{z}) + \mathbf{H}_{1,6}\Delta\mathbf{z}^{[2]} + h.o.t.\big) \otimes \\ \big(\mathbf{H}_{1,1}\Delta\mathbf{x} + \mathbf{H}_{1,2}\Delta\mathbf{x}^{[2]} + \mathbf{H}_{1,4}\Delta\mathbf{z} + \mathbf{H}_{1,5}(\Delta\mathbf{x}\otimes\Delta\mathbf{z}) + \mathbf{H}_{1,6}\Delta\mathbf{z}^{[2]} + h.o.t.\big)$$

$$\begin{aligned}\Delta \boldsymbol{h}^{[2]} = & \big(\mathbf{H}_{1,1}\Delta\mathbf{x}\big)\otimes\big(\mathbf{H}_{1,1}\Delta\mathbf{x}\big) + \big(\mathbf{H}_{1,2}\Delta\mathbf{x}^{[2]}\big)\otimes\big(\mathbf{H}_{1,1}\Delta\mathbf{x}\big) + \big(\mathbf{H}_{1,1}\Delta\mathbf{x}\big)\otimes\big(\mathbf{H}_{1,2}\Delta\mathbf{x}^{[2]}\big) + \\ & + \big(\mathbf{H}_{1,4}\Delta\mathbf{z}\big)\otimes\big(\mathbf{H}_{1,1}\Delta\mathbf{x}\big) + \big(\mathbf{H}_{1,1}\Delta\mathbf{x}\big)\otimes\big(\mathbf{H}_{1,4}\Delta\mathbf{z}\big) + \\ & + \big(\mathbf{H}_{1,5}(\Delta\mathbf{x}\otimes\Delta\mathbf{z})\big)\otimes\big(\mathbf{H}_{1,1}\Delta\mathbf{x}\big) + \big(\mathbf{H}_{1,1}\Delta\mathbf{x}\big)\otimes\big(\mathbf{H}_{1,5}(\Delta\mathbf{x}\otimes\Delta\mathbf{z})\big) + \\ & + \big(\mathbf{H}_{1,6}\Delta\mathbf{z}^{[2]}\big)\otimes\big(\mathbf{H}_{1,1}\Delta\mathbf{x}\big) + \big(\mathbf{H}_{1,1}\Delta\mathbf{x}\big)\otimes\big(\mathbf{H}_{1,6}\Delta\mathbf{z}^{[2]}\big) + \\ & + \big(\mathbf{H}_{1,2}\Delta\mathbf{x}^{[2]}\big)\otimes\big(\mathbf{H}_{1,4}\Delta\mathbf{z}\big) + \big(\mathbf{H}_{1,4}\Delta\mathbf{z}\big)\otimes\big(\mathbf{H}_{1,2}\Delta\mathbf{x}^{[2]}\big) + \big(\mathbf{H}_{1,4}\Delta\mathbf{z}\big)\otimes\big(\mathbf{H}_{1,4}\Delta\mathbf{z}\big) + \\ & + \big(\mathbf{H}_{1,5}(\Delta\mathbf{x}\otimes\Delta\mathbf{z})\big)\otimes\big(\mathbf{H}_{1,4}\Delta\mathbf{z}\big) + \big(\mathbf{H}_{1,4}\Delta\mathbf{z}\big)\otimes\big(\mathbf{H}_{1,5}(\Delta\mathbf{x}\otimes\Delta\mathbf{z})\big) + \\ & + \big(\mathbf{H}_{1,6}\Delta\mathbf{z}^{[2]}\big)\otimes\big(\mathbf{H}_{1,4}\Delta\mathbf{z}\big) + \big(\mathbf{H}_{1,4}\Delta\mathbf{z}\big)\otimes\big(\mathbf{H}_{1,6}\Delta\mathbf{z}^{[2]}\big) + h.o.t.\end{aligned}$$

$$\begin{aligned}\Delta \boldsymbol{h}^{[2]} = & \big(\mathbf{H}_{1,1}\otimes\mathbf{H}_{1,1}\big)\Delta\mathbf{x}^{[2]} + \big(\mathbf{H}_{1,2}\otimes\mathbf{H}_{1,1} + \mathbf{H}_{1,1}\otimes\mathbf{H}_{1,2}\big)\Delta\mathbf{x}^{[3]} + \\ & + \big(\mathbf{H}_{1,4}\otimes\mathbf{H}_{1,1}\big)(\Delta\mathbf{z}\otimes\Delta\mathbf{x}) + \big(\mathbf{H}_{1,1}\otimes\mathbf{H}_{1,4}\big)(\Delta\mathbf{x}\otimes\Delta\mathbf{z}) + \\ & + \big(\mathbf{H}_{1,5}\otimes\mathbf{H}_{1,1}\big)(\Delta\mathbf{x}\otimes\Delta\mathbf{z}\otimes\Delta\mathbf{x}) + \big(\mathbf{H}_{1,1}\otimes\mathbf{H}_{1,5}\big)(\Delta\mathbf{x}^{[2]}\otimes\Delta\mathbf{z}) + \\ & + \big(\mathbf{H}_{1,6}\otimes\mathbf{H}_{1,1}\big)(\Delta\mathbf{z}^{[2]}\otimes\Delta\mathbf{x}) + \big(\mathbf{H}_{1,1}\otimes\mathbf{H}_{1,6}\big)(\Delta\mathbf{x}\otimes\Delta\mathbf{z}^{[2]}) + \\ & + \big(\mathbf{H}_{1,2}\otimes\mathbf{H}_{1,4}\big)(\Delta\mathbf{x}^{[2]}\otimes\Delta\mathbf{z}) + \big(\mathbf{H}_{1,4}\otimes\mathbf{H}_{1,2}\big)(\Delta\mathbf{z}\otimes\Delta\mathbf{x}^{[2]}) + \big(\mathbf{H}_{1,4}\otimes\mathbf{H}_{1,4}\big)\Delta\mathbf{z}^{[2]} + \\ & + \big(\mathbf{H}_{1,5}\otimes\mathbf{H}_{1,4}\big)(\Delta\mathbf{x}\otimes\Delta\mathbf{z}^{[2]}) + \big(\mathbf{H}_{1,4}\otimes\mathbf{H}_{1,5}\big)(\Delta\mathbf{z}\otimes\Delta\mathbf{x}\otimes\Delta\mathbf{z}) + \\ & + \big(\mathbf{H}_{1,6}\otimes\mathbf{H}_{1,4} + \mathbf{H}_{1,4}\otimes\mathbf{H}_{1,6}\big)\Delta\mathbf{z}^{[3]} + h.o.t.\end{aligned} \tag{24}$$

Following equation (24) and in analogy with previous derivations, we introduce the following definitions:

$$\mathbf{H}_{3,2} = \mathbf{H}_{1,1} \otimes \mathbf{H}_{1,1} \tag{25}$$

$$\mathbf{H}_{3,3} = \mathbf{H}_{1,2} \otimes \mathbf{H}_{1,1} + \mathbf{H}_{1,1} \otimes \mathbf{H}_{1,2} \tag{26}$$

$$\mathbf{H}_{3,5}(\Delta\mathbf{x} \otimes \Delta\mathbf{z}) = (\mathbf{H}_{1,4} \otimes \mathbf{H}_{1,1})(\Delta\mathbf{z} \otimes \Delta\mathbf{x}) + (\mathbf{H}_{1,1} \otimes \mathbf{H}_{1,4})(\Delta\mathbf{x} \otimes \Delta\mathbf{z}) \tag{27}$$

$$\mathbf{H}_{3,6} = \mathbf{H}_{1,4} \otimes \mathbf{H}_{1,4} \tag{28}$$

$$\mathbf{H}_{3,7}(\Delta\mathbf{x}^{[2]} \otimes \Delta\mathbf{z}) = (\mathbf{H}_{1,5} \otimes \mathbf{H}_{1,1})(\Delta\mathbf{x} \otimes \Delta\mathbf{z} \otimes \Delta\mathbf{x}) + (\mathbf{H}_{1,1} \otimes \mathbf{H}_{1,5})(\Delta\mathbf{x}^{[2]} \otimes \Delta\mathbf{z}) + \\ + (\mathbf{H}_{1,4} \otimes \mathbf{H}_{1,2})(\Delta\mathbf{z} \otimes \Delta\mathbf{x}^{[2]}) + (\mathbf{H}_{1,2} \otimes \mathbf{H}_{1,4})(\Delta\mathbf{x}^{[2]} \otimes \Delta\mathbf{z}) \tag{29}$$

$$\mathbf{H}_{3,8}(\Delta\mathbf{x} \otimes \Delta\mathbf{z}^{[2]}) = (\mathbf{H}_{1,6} \otimes \mathbf{H}_{1,1})(\Delta\mathbf{z}^{[2]} \otimes \Delta\mathbf{x}) + (\mathbf{H}_{1,1} \otimes \mathbf{H}_{1,6})(\Delta\mathbf{x} \otimes \Delta\mathbf{z}^{[2]}) \\ + (\mathbf{H}_{1,5} \otimes \mathbf{H}_{1,4})(\Delta\mathbf{x} \otimes \Delta\mathbf{z}^{[2]}) + (\mathbf{H}_{1,4} \otimes \mathbf{H}_{1,5})(\Delta\mathbf{z} \otimes \Delta\mathbf{x} \otimes \Delta\mathbf{z}) \tag{30}$$

$$\mathbf{H}_{3,9}\Delta\mathbf{z}^{[3]} = (\mathbf{H}_{1,6} \otimes \mathbf{H}_{1,4} + \mathbf{H}_{1,4} \otimes \mathbf{H}_{1,6})\Delta\mathbf{z}^{[3]} \tag{31}$$

Therefore, utilizing definitions (25) to (31), equation (24) can be reformulated as follows:

$$\mathbf{0} = \Delta h^{[2]} = \mathbf{H}_{3,2}\Delta\mathbf{x}^{[2]} + \mathbf{H}_{3,3}\Delta\mathbf{x}^{[3]} + \mathbf{H}_{3,5}(\Delta\mathbf{x} \otimes \Delta\mathbf{z}) + \\ + \mathbf{H}_{3,6}\Delta\mathbf{z}^{[2]} + \mathbf{H}_{3,7}(\Delta\mathbf{x}^{[2]} \otimes \Delta\mathbf{z}) + \mathbf{H}_{3,8}(\Delta\mathbf{x} \otimes \Delta\mathbf{z}^{[2]}) + \mathbf{H}_{3,9}\Delta\mathbf{z}^{[3]} + h.o.t. \tag{32}$$

Based on equation (32), the term $\Delta h^{[3]}$ is derived as follows:

$$\mathbf{0} = \Delta h^{[3]} = \Delta h^{[2]} \otimes \Delta h = (\mathbf{H}_{3,2}\Delta\mathbf{x}^{[2]} + \mathbf{H}_{3,3}\Delta\mathbf{x}^{[3]} + \mathbf{H}_{3,5}(\Delta\mathbf{x} \otimes \Delta\mathbf{z}) + \\ + \mathbf{H}_{3,6}\Delta\mathbf{z}^{[2]} + \mathbf{H}_{3,7}(\Delta\mathbf{x}^{[2]} \otimes \Delta\mathbf{z}) + \mathbf{H}_{3,8}(\Delta\mathbf{x} \otimes \Delta\mathbf{z}^{[2]}) + \mathbf{H}_{3,9}\Delta\mathbf{z}^{[3]} + h.o.t.) \otimes \\ (\mathbf{H}_{1,1}\Delta\mathbf{x} + \mathbf{H}_{1,2}\Delta\mathbf{x}^{[2]} + \mathbf{H}_{1,4}\Delta\mathbf{z} + \mathbf{H}_{1,5}(\Delta\mathbf{x} \otimes \Delta\mathbf{z}) + \mathbf{H}_{1,6}\Delta\mathbf{z}^{[2]} + h.o.t.)$$

$$\Delta h^{[3]} = (\mathbf{H}_{3,2}\Delta\mathbf{x}^{[2]}) \otimes (\mathbf{H}_{1,1}\Delta\mathbf{x}) + (\mathbf{H}_{3,2}\Delta\mathbf{x}^{[2]}) \otimes (\mathbf{H}_{1,4}\Delta\mathbf{z}) + \\ + (\mathbf{H}_{3,5}(\Delta\mathbf{x} \otimes \Delta\mathbf{z})) \otimes (\mathbf{H}_{1,1}\Delta\mathbf{x}) + (\mathbf{H}_{3,5}(\Delta\mathbf{x} \otimes \Delta\mathbf{z})) \otimes (\mathbf{H}_{1,4}\Delta\mathbf{z}) + \\ + (\mathbf{H}_{3,6}\Delta\mathbf{z}^{[2]}) \otimes (\mathbf{H}_{1,1}\Delta\mathbf{x}) + (\mathbf{H}_{3,6}\Delta\mathbf{z}^{[2]}) \otimes (\mathbf{H}_{1,4}\Delta\mathbf{z}) + h.o.t.$$

$$\Delta h^{[3]} = (\mathbf{H}_{3,2} \otimes \mathbf{H}_{1,1})\Delta\mathbf{x}^{[3]} + (\mathbf{H}_{3,2} \otimes \mathbf{H}_{1,4})(\Delta\mathbf{x}^{[2]} \otimes \Delta\mathbf{z}) + \\ + (\mathbf{H}_{3,5} \otimes \mathbf{H}_{1,1})(\Delta\mathbf{x} \otimes \Delta\mathbf{z} \otimes \Delta\mathbf{x}) + (\mathbf{H}_{3,5} \otimes \mathbf{H}_{1,4})(\Delta\mathbf{x} \otimes \Delta\mathbf{z}^{[2]}) + \\ + (\mathbf{H}_{3,6} \otimes \mathbf{H}_{1,1})(\Delta\mathbf{z}^{[2]} \otimes \Delta\mathbf{x}) + (\mathbf{H}_{3,6} \otimes \mathbf{H}_{1,4})\Delta\mathbf{z}^{[3]} + h.o.t. \tag{33}$$

where the following terms can be defined:

$$\mathbf{H}_{6,3} = \mathbf{H}_{3,2} \otimes \mathbf{H}_{1,1} \tag{34}$$

$$\mathbf{H}_{6,7}(\Delta\mathbf{x}^{[2]} \otimes \Delta\mathbf{z}) = (\mathbf{H}_{3,2} \otimes \mathbf{H}_{1,4})(\Delta\mathbf{x}^{[2]} \otimes \Delta\mathbf{z}) + (\mathbf{H}_{3,5} \otimes \mathbf{H}_{1,1})(\Delta\mathbf{x} \otimes \Delta\mathbf{z} \otimes \Delta\mathbf{x}) \quad (35)$$

$$\mathbf{H}_{6,8}(\Delta\mathbf{x} \otimes \Delta\mathbf{z}^{[2]}) = (\mathbf{H}_{3,5} \otimes \mathbf{H}_{1,4})(\Delta\mathbf{x} \otimes \Delta\mathbf{z}^{[2]}) + (\mathbf{H}_{3,6} \otimes \mathbf{H}_{1,1})(\Delta\mathbf{z}^{[2]} \otimes \Delta\mathbf{x}) \quad (36)$$

$$\mathbf{H}_{6,9} = \mathbf{H}_{3,6} \otimes \mathbf{H}_{1,4} \quad (37)$$

Then, equation (33) can be rewritten as follows:

$$\begin{aligned}\mathbf{0} = \Delta\boldsymbol{h}^{[3]} &= \mathbf{H}_{6,3}\Delta\mathbf{x}^{[3]} + \mathbf{H}_{6,7}(\Delta\mathbf{x}^{[2]} \otimes \Delta\mathbf{z}) + \\ &+ \mathbf{H}_{6,8}(\Delta\mathbf{x} \otimes \Delta\mathbf{z}^{[2]}) + \mathbf{H}_{6,9}\Delta\mathbf{z}^{[3]} + h.o.t.\end{aligned} \quad (38)$$

On the other hand, for the term $\Delta\boldsymbol{h} \otimes \Delta\mathbf{x}$ we derive the following sequence:

$$\begin{aligned}\mathbf{0} = \Delta\boldsymbol{h} \otimes \Delta\mathbf{x} = \big(&\mathbf{H}_{1,1}\Delta\mathbf{x} + \mathbf{H}_{1,2}\Delta\mathbf{x}^{[2]} + \mathbf{H}_{1,3}\Delta\mathbf{x}^{[3]} + \mathbf{H}_{1,4}\Delta\mathbf{z} + \mathbf{H}_{1,5}(\Delta\mathbf{x} \otimes \Delta\mathbf{z}) + \\ &+ \mathbf{H}_{1,6}\Delta\mathbf{z}^{[2]} + \mathbf{H}_{1,7}(\Delta\mathbf{x}^{[2]} \otimes \Delta\mathbf{z}) + \mathbf{H}_{1,8}(\Delta\mathbf{x} \otimes \Delta\mathbf{z}^{[2]}) + \mathbf{H}_{1,9}\Delta\mathbf{z}^{[3]} + h.o.t.\big) \otimes \Delta\mathbf{x}\end{aligned}$$

$$\begin{aligned}\mathbf{0} = \Delta\boldsymbol{h} \otimes \Delta\mathbf{x} &= (\mathbf{H}_{1,1}\Delta\mathbf{x}) \otimes \Delta\mathbf{x} + (\mathbf{H}_{1,2}\Delta\mathbf{x}^{[2]}) \otimes \Delta\mathbf{x} + (\mathbf{H}_{1,4}\Delta\mathbf{z}) \otimes \Delta\mathbf{x} + \\ &+ \big(\mathbf{H}_{1,5}(\Delta\mathbf{x} \otimes \Delta\mathbf{z})\big) \otimes \Delta\mathbf{x} + (\mathbf{H}_{1,6}\Delta\mathbf{z}^{[2]}) \otimes \Delta\mathbf{x} + h.o.t.\end{aligned}$$

$$\begin{aligned}\mathbf{0} = \Delta\boldsymbol{h} \otimes \Delta\mathbf{x} &= (\mathbf{H}_{1,1} \otimes \mathbf{I}_N)\Delta\mathbf{x}^{[2]} + (\mathbf{H}_{1,2} \otimes \mathbf{I}_N)\Delta\mathbf{x}^{[3]} + (\mathbf{H}_{1,4} \otimes \mathbf{I}_N)(\Delta\mathbf{z} \otimes \Delta\mathbf{x}) + \\ &+ (\mathbf{H}_{1,5} \otimes \mathbf{I}_N)(\Delta\mathbf{x} \otimes \Delta\mathbf{z} \otimes \Delta\mathbf{x}) + (\mathbf{H}_{1,6} \otimes \mathbf{I}_N)(\Delta\mathbf{z}^{[2]} \otimes \Delta\mathbf{x}) + h.o.t.\end{aligned} \quad (39)$$

Then, by defining

$$\mathbf{H}_{2,2} = \mathbf{H}_{1,1} \otimes \mathbf{I}_N \quad (40)$$

$$\mathbf{H}_{2,3} = \mathbf{H}_{1,2} \otimes \mathbf{I}_N \quad (41)$$

$$\mathbf{H}_{2,5}(\Delta\mathbf{x} \otimes \Delta\mathbf{z}) = (\mathbf{H}_{1,4} \otimes \mathbf{I}_N)(\Delta\mathbf{z} \otimes \Delta\mathbf{x}) \quad (42)$$

$$\mathbf{H}_{2,7}(\Delta\mathbf{x}^{[2]} \otimes \Delta\mathbf{z}) = (\mathbf{H}_{1,5} \otimes \mathbf{I}_N)(\Delta\mathbf{x} \otimes \Delta\mathbf{z} \otimes \Delta\mathbf{x}) \quad (43)$$

$$\mathbf{H}_{2,8}(\Delta\mathbf{x} \otimes \Delta\mathbf{z}^{[2]}) = (\mathbf{H}_{1,6} \otimes \mathbf{I}_N)(\Delta\mathbf{z}^{[2]} \otimes \Delta\mathbf{x}) \quad (44)$$

Equation (39) can be re-expressed as follows:

$$\begin{aligned}\mathbf{0} = \Delta\boldsymbol{h} \otimes \Delta\mathbf{x} &= \mathbf{H}_{2,2}\Delta\mathbf{x}^{[2]} + \mathbf{H}_{2,3}\Delta\mathbf{x}^{[3]} + \mathbf{H}_{2,5}(\Delta\mathbf{x} \otimes \Delta\mathbf{z}) + \\ &+ \mathbf{H}_{2,7}(\Delta\mathbf{x}^{[2]} \otimes \Delta\mathbf{z}) + \mathbf{H}_{2,8}(\Delta\mathbf{x} \otimes \Delta\mathbf{z}^{[2]}) + h.o.t.\end{aligned} \quad (45)$$

Moreover, for the term $\Delta\boldsymbol{h} \otimes \Delta\mathbf{x}^{[2]}$, we have

$$\begin{aligned}\mathbf{0} = \Delta\boldsymbol{h} \otimes \Delta\mathbf{x}^{[2]} = \big(&\mathbf{H}_{1,1}\Delta\mathbf{x} + \mathbf{H}_{1,2}\Delta\mathbf{x}^{[2]} + \mathbf{H}_{1,3}\Delta\mathbf{x}^{[3]} + \mathbf{H}_{1,4}\Delta\mathbf{z} + \mathbf{H}_{1,5}(\Delta\mathbf{x} \otimes \Delta\mathbf{z}) + \\ &+ \mathbf{H}_{1,6}\Delta\mathbf{z}^{[2]} + \mathbf{H}_{1,7}(\Delta\mathbf{x}^{[2]} \otimes \Delta\mathbf{z}) + \mathbf{H}_{1,8}(\Delta\mathbf{x} \otimes \Delta\mathbf{z}^{[2]}) + \mathbf{H}_{1,9}\Delta\mathbf{z}^{[3]} + h.o.t.\big) \otimes \Delta\mathbf{x}^{[2]}\end{aligned}$$

$$\mathbf{0} = \Delta\boldsymbol{h} \otimes \Delta\mathbf{x}^{[2]} = (\mathbf{H}_{1,1}\Delta\mathbf{x}) \otimes \Delta\mathbf{x}^{[2]} + (\mathbf{H}_{1,4}\Delta\mathbf{z}) \otimes \Delta\mathbf{x}^{[2]} + h.o.t.$$

$$\mathbf{0} = \Delta \boldsymbol{h} \otimes \Delta \mathbf{x}^{[2]} = \left(\mathbf{H}_{1,1} \otimes \mathbf{I}_{N^2}\right) \Delta \mathbf{x}^{[3]} + \left(\mathbf{H}_{1,4} \otimes \mathbf{I}_{N^2}\right)\left(\Delta \mathbf{z} \otimes \Delta \mathbf{x}^{[2]}\right) + h.o.t. \quad (46)$$

In equation (46), by defining:

$$\mathbf{H}_{4,3} = \mathbf{H}_{1,1} \otimes \mathbf{I}_{N^2} \quad (47)$$

$$\mathbf{H}_{4,7}\left(\Delta \mathbf{x}^{[2]} \otimes \Delta \mathbf{z}\right) = \left(\mathbf{H}_{1,4} \otimes \mathbf{I}_{N^2}\right)\left(\Delta \mathbf{z} \otimes \Delta \mathbf{x}^{[2]}\right) \quad (48)$$

it follows that:

$$\mathbf{0} = \Delta \boldsymbol{h} \otimes \Delta \mathbf{x}^{[2]} = \mathbf{H}_{4,3} \Delta \mathbf{x}^{[3]} + \mathbf{H}_{4,7}\left(\Delta \mathbf{x}^{[2]} \otimes \Delta \mathbf{z}\right) + h.o.t. \quad (49)$$

Finally, for the term $\Delta \boldsymbol{h}^{[2]} \otimes \Delta \mathbf{x}$, it can be found that:

$$\mathbf{0} = \Delta \boldsymbol{h}^{[2]} \otimes \Delta \mathbf{x} = \left(\mathbf{H}_{3,2} \Delta \mathbf{x}^{[2]} + \mathbf{H}_{3,3} \Delta \mathbf{x}^{[3]} + \mathbf{H}_{3,5}(\Delta \mathbf{x} \otimes \Delta \mathbf{z}) + \mathbf{H}_{3,6} \Delta \mathbf{z}^{[2]} + \right.$$
$$\left. + \mathbf{H}_{3,7}\left(\Delta \mathbf{x}^{[2]} \otimes \Delta \mathbf{z}\right) + \mathbf{H}_{3,8}\left(\Delta \mathbf{x} \otimes \Delta \mathbf{z}^{[2]}\right) + \mathbf{H}_{3,9} \Delta \mathbf{z}^{[3]} + h.o.t. \right) \otimes \Delta \mathbf{x}$$

$$\mathbf{0} = \Delta \boldsymbol{h}^{[2]} \otimes \Delta \mathbf{x} = \left(\mathbf{H}_{3,2} \Delta \mathbf{x}^{[2]}\right) \otimes \Delta \mathbf{x} + \left(\mathbf{H}_{3,5}(\Delta \mathbf{x} \otimes \Delta \mathbf{z})\right) \otimes \Delta \mathbf{x} +$$
$$+ \left(\mathbf{H}_{3,6} \Delta \mathbf{z}^{[2]}\right) \otimes \Delta \mathbf{x} + h.o.t.$$

$$\mathbf{0} = \Delta \boldsymbol{h}^{[2]} \otimes \Delta \mathbf{x} = \left(\mathbf{H}_{3,2} \otimes \mathbf{I}_N\right) \Delta \mathbf{x}^{[3]} + \left(\mathbf{H}_{3,5} \otimes \mathbf{I}_N\right)(\Delta \mathbf{x} \otimes \Delta \mathbf{z} \otimes \Delta \mathbf{x}) +$$
$$+ \left(\mathbf{H}_{3,6} \otimes \mathbf{I}_N\right)\left(\Delta \mathbf{z}^{[2]} \otimes \Delta \mathbf{x}\right) + h.o.t. \quad (50)$$

In (50), the following definitions can be introduced:

$$\mathbf{H}_{5,3} = \mathbf{H}_{3,2} \otimes \mathbf{I}_N \quad (51)$$

$$\mathbf{H}_{5,7}\left(\Delta \mathbf{z} \otimes \Delta \mathbf{x}^{[2]}\right) = \left(\mathbf{H}_{3,5} \otimes \mathbf{I}_N\right)(\Delta \mathbf{x} \otimes \Delta \mathbf{z} \otimes \Delta \mathbf{x}) \quad (52)$$

$$\mathbf{H}_{5,8}\left(\Delta \mathbf{x} \otimes \Delta \mathbf{z}^{[2]}\right) = \left(\mathbf{H}_{3,6} \otimes \mathbf{I}_N\right)\left(\Delta \mathbf{z}^{[2]} \otimes \Delta \mathbf{x}\right) \quad (53)$$

Substituting definitions (51)-(53) into equation (50), the expression $\Delta \boldsymbol{h}^{[2]} \otimes \Delta \mathbf{x}$ becomes:

$$\mathbf{0} = \Delta \boldsymbol{h}^{[2]} \otimes \Delta \mathbf{x} = \mathbf{H}_{5,3} \Delta \mathbf{x}^{[3]} + \mathbf{H}_{5,7}\left(\Delta \mathbf{z} \otimes \Delta \mathbf{x}^{[2]}\right) +$$
$$+ \mathbf{H}_{5,8}\left(\Delta \mathbf{x} \otimes \Delta \mathbf{z}^{[2]}\right) + h.o.t. \quad (54)$$

Now, by incorporating (32), (38), (45), (49), and (54) into (19), we obtain:

$$\begin{bmatrix} \Delta \dot{\mathbf{x}} \\ \Delta \dot{\mathbf{x}}^{[2]} \\ \Delta \dot{\mathbf{x}}^{[3]} \\ 0 \\ 0 \\ 0 \\ 0 \\ 0 \\ 0 \\ \vdots \end{bmatrix} = \begin{bmatrix} G_{1,1} & G_{1,2} & G_{1,3} & G_{1,4} & G_{1,5} & G_{1,6} & G_{1,7} & G_{1,8} & G_{1,9} \\ 0 & G_{2,2} & G_{2,3} & 0 & G_{2,5} & 0 & G_{2,7} & G_{2,8} & 0 \\ 0 & 0 & G_{3,3} & 0 & 0 & 0 & G_{3,7} & 0 & 0 \\ H_{1,1} & H_{1,2} & H_{1,3} & H_{1,4} & H_{1,5} & H_{1,6} & H_{1,7} & H_{1,8} & H_{1,9} \\ 0 & H_{2,2} & H_{2,3} & 0 & H_{2,5} & 0 & H_{2,7} & H_{2,8} & 0 \\ 0 & H_{3,2} & H_{3,3} & 0 & H_{3,5} & H_{3,6} & H_{3,7} & H_{3,8} & H_{3,9} \\ 0 & 0 & H_{4,3} & 0 & 0 & 0 & H_{4,7} & 0 & 0 \\ 0 & 0 & H_{5,3} & 0 & 0 & 0 & H_{5,7} & H_{5,8} & 0 \\ 0 & 0 & H_{6,3} & 0 & 0 & 0 & H_{6,7} & H_{6,8} & H_{6,9} \\ \vdots & \vdots & \vdots & \vdots & \vdots & \vdots & \vdots & \vdots & \vdots \end{bmatrix} \begin{bmatrix} \Delta \mathbf{x} \\ \Delta \mathbf{x}^{[2]} \\ \Delta \mathbf{x}^{[3]} \\ \Delta \mathbf{z} \\ \Delta \mathbf{x} \otimes \Delta \mathbf{z} \\ \Delta \mathbf{z}^{[2]} \\ \Delta \mathbf{x}^{[2]} \otimes \Delta \mathbf{z} \\ \Delta \mathbf{x} \otimes \Delta \mathbf{z}^{[2]} \\ \Delta \mathbf{z}^{[3]} \end{bmatrix} \quad (55)$$

$+h.o.t.$

The system in (55) can be further simplified by neglecting terms of order four and higher, resulting in:

$$\begin{bmatrix} \Delta \dot{\mathbf{x}}_{3ord} \\ 0 \end{bmatrix} = \begin{bmatrix} \mathbf{F}_{11} & \mathbf{F}_{12} \\ \mathbf{F}_{21} & \mathbf{F}_{22} \end{bmatrix} \begin{bmatrix} \Delta \mathbf{x}_{3ord} \\ \Delta \mathbf{z}_{3ord} \end{bmatrix} \quad (56)$$

where the vectors $\Delta \mathbf{x}_{3ord}$ and $\Delta \mathbf{z}_{3ord}$ are defined as:

$$\Delta \mathbf{x}_{3ord} = \begin{bmatrix} \Delta \mathbf{x} \\ \Delta \mathbf{x}^{[2]} \\ \Delta \mathbf{x}^{[3]} \end{bmatrix}, \quad \Delta \mathbf{z}_{3ord} = \begin{bmatrix} \Delta \mathbf{z} \\ \Delta \mathbf{x} \otimes \Delta \mathbf{z} \\ \Delta \mathbf{z}^{[2]} \\ \Delta \mathbf{x}^{[2]} \otimes \Delta \mathbf{z} \\ \Delta \mathbf{x} \otimes \Delta \mathbf{z}^{[2]} \\ \Delta \mathbf{z}^{[3]} \end{bmatrix}, \quad (57)$$

The submatrices $\mathbf{F}_{11}$, $\mathbf{F}_{12}$, $\mathbf{F}_{21}$, and $\mathbf{F}_{22}$ are structured as follows:

$$\mathbf{F}_{11} = \begin{bmatrix} G_{1,1} & G_{1,3} & G_{1,2} \\ 0 & G_{2,2} & G_{2,3} \\ 0 & 0 & G_{3,3} \end{bmatrix} \quad (58)$$

$$\mathbf{F}_{12} = \begin{bmatrix} G_{1,4} & G_{1,5} & G_{1,6} & G_{1,7} & G_{1,8} & G_{1,9} \\ 0 & G_{2,5} & 0 & G_{2,7} & G_{2,8} & G_{2,9} \\ 0 & 0 & 0 & G_{3,7} & 0 & G_{3,3} \end{bmatrix} \quad (59)$$

$$\mathbf{F}_{21} = \begin{bmatrix} H_{1,1} & H_{1,2} & H_{1,3} \\ 0 & H_{2,2} & H_{2,3} \\ 0 & H_{3,2} & H_{3,3} \\ 0 & 0 & H_{4,3} \\ 0 & 0 & H_{5,3} \\ 0 & 0 & H_{6,3} \end{bmatrix} \quad (60)$$

$$\mathbf{F}_{22} = \begin{bmatrix} \mathbf{H}_{1,4} & \mathbf{H}_{1,5} & \mathbf{H}_{1,6} & \mathbf{H}_{1,7} & \mathbf{H}_{1,8} & \mathbf{H}_{1,9} \\ 0 & \mathbf{H}_{2,5} & 0 & \mathbf{H}_{2,7} & \mathbf{H}_{2,8} & 0 \\ 0 & \mathbf{H}_{3,5} & \mathbf{H}_{3,6} & \mathbf{H}_{3,7} & \mathbf{H}_{3,8} & \mathbf{H}_{3,9} \\ 0 & 0 & 0 & \mathbf{H}_{4,7} & 0 & 0 \\ 0 & 0 & 0 & \mathbf{H}_{5,7} & \mathbf{H}_{5,8} & 0 \\ 0 & 0 & 0 & \mathbf{H}_{6,7} & \mathbf{H}_{6,8} & \mathbf{H}_{6,9} \end{bmatrix} \quad (61)$$

Furthermore, if a Kron reduction is applied to (56) in the following manner:

$$\widetilde{\mathbf{F}}_{11} = \mathbf{F}_{11} - \mathbf{F}_{12}(\mathbf{F}_{22})^{-1}\mathbf{F}_{21}, \quad (62)$$

then an equivalent ODE representation can be expressed as:

$$\Delta \dot{\mathbf{x}}_{3ord} = \widetilde{\mathbf{F}}_{11} \Delta \mathbf{x}_{3ord} \quad (63)$$

It must be emphasized that, from (56) and (62), the matrices $\widetilde{\mathbf{H}}_{1,j}$ (i.e., the reduced-order approximations of the nonlinear terms) can be obtained from the first $M$ rows of the matrix $(\mathbf{F}_{22})^{-1}\mathbf{F}_{21}$.

As the reader may observe, a critical condition for achieving the equivalent ODE form (63) is the invertibility of the matrix $\mathbf{F}_{22}$. This issue is addressed in the following Theorem.

Theorem 2. Consider a nonlinear DAE system of the form (8), (9) that has been linearized up to the third order as in equation (20). If the functions $\Delta \boldsymbol{h}^{[2]}, \Delta \boldsymbol{h}^{[3]}, \Delta \boldsymbol{h} \otimes \Delta \mathbf{x}, \Delta \boldsymbol{h}^{[2]} \otimes \Delta \mathbf{x}$, and $\Delta \boldsymbol{h} \otimes \Delta \mathbf{x}^{[2]}$ are incorporated for obtaining a third-order representation of the form (56), then the matrix $\mathbf{F}_{22}$ is invertible if and only if $\det(\mathbf{H}_{1,4}) \neq 0$.

*Proof*:

According to [26], the matrix $\mathbf{F}_{22}$ is invertible if and only if $\det(\mathbf{F}_{22}) \neq 0$. The matrix $\mathbf{F}_{22}$ can then be rewritten as follows:

$$\mathbf{F}_{22} = \begin{bmatrix} \widehat{\mathbf{F}}_{11} & \widehat{\mathbf{F}}_{12} & \widehat{\mathbf{F}}_{13} \\ 0 & \widehat{\mathbf{F}}_{22} & \widehat{\mathbf{F}}_{23} \\ 0 & 0 & \widehat{\mathbf{F}}_{33} \end{bmatrix} \quad (64)$$

where

$$\widehat{\mathbf{F}}_{11} = \mathbf{H}_{1,4}, \quad \widehat{\mathbf{F}}_{12} = [\mathbf{H}_{1,5} \quad \mathbf{H}_{1,6}], \quad \widehat{\mathbf{F}}_{13} = [\mathbf{H}_{1,7} \quad \mathbf{H}_{1,8} \quad \mathbf{H}_{1,9}]$$

$$\widehat{\mathbf{F}}_{22} = \begin{bmatrix} \mathbf{H}_{2,5} & \mathbf{0} \\ \mathbf{H}_{3,5} & \mathbf{H}_{3,6} \end{bmatrix}, \quad \widehat{\mathbf{F}}_{23} = \begin{bmatrix} \mathbf{H}_{2,7} & \mathbf{H}_{2,8} & \mathbf{0} \\ \mathbf{H}_{3,7} & \mathbf{H}_{3,8} & \mathbf{H}_{3,9} \end{bmatrix}, \quad \widehat{\mathbf{F}}_{33} = \begin{bmatrix} \mathbf{H}_{4,7} & \mathbf{0} & \mathbf{0} \\ \mathbf{H}_{5,7} & \mathbf{H}_{5,8} & \mathbf{0} \\ \mathbf{H}_{6,7} & \mathbf{H}_{6,8} & \mathbf{H}_{6,9} \end{bmatrix}$$

As can be noted, the matrix $\mathbf{F}_{22}$ is an upper-triangular block-matrix. According to [26], the determinant of $\mathbf{F}_{22}$ is:

$$\det(\mathbf{F}_{22}) = \det(\mathbf{H}_{1,4}) \det(\widehat{\mathbf{F}}_{22}) \det(\widehat{\mathbf{F}}_{33}) \tag{65}$$

Then, in the same manner we have that

$$\det(\widehat{\mathbf{F}}_{22}) = \det(\mathbf{H}_{2,5}) \det(\mathbf{H}_{3,6}) \tag{66}$$

$$\det(\widehat{\mathbf{F}}_{33}) = \det(\mathbf{H}_{4,7}) \det(\mathbf{H}_{5,8}) \det(\mathbf{H}_{6,9}) \tag{67}$$

Now, by applying the definitions of matrices $\mathbf{H}_{2,5}$, $\mathbf{H}_{3,6}$, $\mathbf{H}_{4,7}$, $\mathbf{H}_{5,8}$, and $\mathbf{H}_{6,9}$ in equations (42), (28), (48), (53), and (37), respectively, and by invoking the properties of the determinant and the Kronecker product, it can be readily shown that:

$$\det(\mathbf{H}_{2,5}) = \det(\mathbf{H}_{1,4})^N \tag{68}$$

$$\det(\mathbf{H}_{3,6}) = \det(\mathbf{H}_{1,4})^{2M} \tag{69}$$

$$\det(\mathbf{H}_{4,7}) = \det(\mathbf{H}_{1,4})^{N^2} \tag{70}$$

$$\det(\mathbf{H}_{5,8}) = \det(\mathbf{H}_{1,4})^{2MN} \tag{71}$$

$$\det(\mathbf{H}_{6,9}) = \det(\mathbf{H}_{3,6})^M \det(\mathbf{H}_{1,4})^{M^2} \tag{72}$$

As can be observed, the determinants in equations (65) through (72) depend on the determinant of matrix $\mathbf{H}_{1,4}$. Therefore, if $\det(\mathbf{H}_{1,4}) \neq 0$, it follows that $\det(\mathbf{F}_{22}) \neq 0$ and hence the matrix $\mathbf{F}_{22}$ is invertible.

∎

In practice, although $\mathbf{F}_{22}$ is theoretically invertible according to Theorem 2, the numerical stability of the inversion depends on the magnitude of $|\det(\mathbf{H}_{1,4})|$. If this determinant is either too small or too large, the value of $\det(\mathbf{F}_{22})$ may also become correspondingly small or large, potentially rendering matrix $\mathbf{F}_{22}$ ill-conditioned. This can lead to significant inaccuracies when computing its inverse [27].

In the following section, numerical examples are presented to illustrate the application of the proposed methodology.

## 4. Numerical Application

*4.1 First Test System*

We begin with a simple and compact test system, defined by the following set of expressions:

$$\begin{aligned}\dot{x}_1 &= 0.2x_1^2 - 2x_1 - 0.5x_2 + 0.1z \\ \dot{x}_2 &= -0.15x_2^2 + 0.1x_1x_2 + 0.3x_1^2 - 2.2x_2 - 2x_1 - 0.1z\end{aligned} \quad (73)$$

$$z = x_1 - x_2 \quad (74)$$

As can be observed, in equation (74), the function $\widetilde{h}$ is already defined and can be directly substituted into equation (73), yielding:

$$\begin{aligned}\dot{x}_1 &= 0.2x_1^2 - 1.9x_1 - 0.6x_2 \\ \dot{x}_2 &= -0.15x_2^2 + 0.1x_1x_2 + 0.3x_1^2 - 2.1x_2 - 2.1x_1\end{aligned} \quad (75)$$

Upon examining equation (75), the reader will notice that this first test system is of order 2, and therefore, the Carleman linearization will result in a system of the form given by equation (6), up to the second-order approximation.

To begin this process, the first step is to determine the equilibrium point of the system, which is given by $\mathbf{x}_{sep} = [0 \quad 0]^T$, $\mathbf{z}_{sep} = 0$.

The resulting second-order Carleman extended model, in the form of equation (6), yields the following linearized system:

$$\begin{bmatrix} \Delta\dot{x}_1 \\ \Delta\dot{x}_2 \\ \Delta\dot{x}_{1,1} \\ \Delta\dot{x}_{1,2} \\ \Delta\dot{x}_{2,1} \\ \Delta\dot{x}_{2,2} \end{bmatrix} = \begin{bmatrix} -1.9 & -0.6 & 0.2 & 0 & 0 & 0 \\ -2.1 & -2.1 & 0.3 & 0.05 & 0.05 & -0.15 \\ 0 & 0 & -3.8 & -0.6 & -0.6 & 0 \\ 0 & 0 & -2.1 & -4 & 0 & -0.6 \\ 0 & 0 & -2.1 & 0 & -4 & -0.6 \\ 0 & 0 & 0 & -2.1 & -2.1 & -4.2 \end{bmatrix} \begin{bmatrix} \Delta x_1 \\ \Delta x_2 \\ \Delta x_1^2 \\ \Delta x_1 \Delta x_2 \\ \Delta x_2 \Delta x_1 \\ \Delta x_2^2 \end{bmatrix} \quad (76)$$

Now, for the DAE system described by equations (73) and (74), if a second-order Carleman-type linearized system of the form (55) is constructed, the result is:

$$[\Delta\dot{x}_1 \ \Delta\dot{x}_2 \ \Delta\dot{x}_{1,1} \ \Delta\dot{x}_{1,2} \ \Delta\dot{x}_{2,1} \ \Delta\dot{x}_{2,2} \ 0 \ 0 \ 0 \ 0]^T =$$

$$\begin{bmatrix} -2 & -0.5 & 0.2 & 0 & 0 & 0 & 0.1 & 0 & 0 & 0 \\ -2 & -2.2 & 0.3 & 0.05 & 0.05 & -0.15 & -0.1 & 0 & 0 & 0 \\ 0 & 0 & -4 & -0.5 & -0.5 & 0 & 0 & 0.2 & 0 & 0 \\ 0 & 0 & -2 & -4.2 & 0 & -0.5 & 0 & -0.1 & 0.1 & 0 \\ 0 & 0 & -2 & 0 & -4.2 & -0.5 & 0 & -0.1 & 0.1 & 0 \\ 0 & 0 & 0 & -2 & -2 & -4.4 & 0 & 0 & -0.2 & 0 \\ -1 & 1 & 0 & 0 & 0 & 0 & 1 & 0 & 0 & 0 \\ 0 & 0 & -1 & 0 & 1 & 0 & 0 & 1 & 0 & 0 \\ 0 & 0 & 0 & -1 & 0 & 1 & 0 & 0 & 1 & 0 \\ 0 & 0 & 1 & -1 & -1 & 1 & 0 & -2 & 2 & 1 \end{bmatrix} \begin{bmatrix} \Delta x_1 \\ \Delta x_2 \\ \Delta x_1^2 \\ \Delta x_1 \Delta x_2 \\ \Delta x_2 \Delta x_1 \\ \Delta x_2^2 \\ \Delta z \\ \Delta x_1 \Delta z \\ \Delta x_2 \Delta z \\ \Delta z^2 \end{bmatrix} \quad (77)$$

In constructing system (77), the second-order matrices $\mathbf{H}_{2,j}$ and $\mathbf{H}_{3,j}$ were directly obtained using equations (25), (27), (28), (40), and (42). Applying the Kron reduction procedure described in equation (62) to equation (77), the resulting equivalent ODE system is:

$$\begin{bmatrix} \Delta\dot{x}_1 \\ \Delta\dot{x}_2 \\ \Delta\dot{x}_{1,1} \\ \Delta\dot{x}_{1,2} \\ \Delta\dot{x}_{2,1} \\ \Delta\dot{x}_{2,2} \end{bmatrix} = \begin{bmatrix} -1.9 & -0.6 & 0.2 & 0 & 0 & 0 \\ -2.1 & -2.1 & 0.3 & 0.05 & 0.05 & -0.15 \\ 0 & 0 & -3.8 & -0.5 & -0.7 & 0 \\ 0 & 0 & -2.1 & -4.1 & 0.1 & -0.6 \\ 0 & 0 & -2.1 & 0.1 & -4.1 & -0.6 \\ 0 & 0 & 0 & -2.2 & -2.0 & -4.2 \end{bmatrix} \begin{bmatrix} \Delta x_1 \\ \Delta x_2 \\ \Delta x_1^2 \\ \Delta x_1 \Delta x_2 \\ \Delta x_2 \Delta x_1 \\ \Delta x_2^2 \end{bmatrix} \quad (78)$$

By comparing equations (76) and (78), it can be observed that there are slight differences in the coefficients associated with the terms $\Delta x_1 \Delta x_2$ and $\Delta x_2 \Delta x_1$. However, since these two terms are mathematically equivalent, the systems described by equations (76) and (78) are functionally identical.

*4.2 Second Test System*

As a second test system, we consider the DAE representation defined by equations (79) and (80) below, where the system consists of three state variables and two algebraic variables:

$$\begin{aligned} \dot{x}_1 &= \frac{x_2}{x_1+1} + e^{-z_1} - x_3^3 + \sin(x_1 z_2) \\ \dot{x}_2 &= \tan(x_1 + z_2) - x_2^2 + x_3 e^{x_1} - \frac{x_1}{x_3+1} \\ \dot{x}_3 &= x_1 x_2 - x_3 + \frac{z_1}{z_1^2+1} + \cos(x_3) \end{aligned} \quad (79)$$

$$\begin{aligned} 0 &= z_1 + x_1 - x_3 \\ 0 &= z_2 + x_2 \end{aligned} \quad (80)$$

By substituting equations (80) into (79), an ODE-only representation is obtained, as given in equation (81):

$$\begin{aligned} \dot{x}_1 &= \frac{x_2}{x_1+1} + e^{x_1-x_3} - x_3^3 + \sin(x_1(1-x_2)) \\ \dot{x}_2 &= \tan(x_1 - x_2 + 1) - x_2^2 + x_3 e^{x_1} - \frac{x_1}{x_3+1} \\ \dot{x}_3 &= x_1 x_2 - x_3 + \frac{x_3-x_1}{(x_3-x_1)^2+1} + \cos(x_3) \end{aligned} \quad (81)$$

In this case, the equilibrium point the system is $\mathbf{x}_{sep} = [0.1369 \quad 1.10817 \quad 1.102644]^T$ and $\mathbf{z}_{sep} = [0.965743 \quad -0.10817]^T$.

For the system described by equations (79) and (80), a third-order equivalent ODE representation will be computed. To illustrate the derivation of coefficient matrices within this representation, let us consider the matrix $\mathbf{G}_{1,4}$, which in this case is given by:

$$\mathbf{G}_{1,4} = \begin{bmatrix} -0.3807 & 0.1369 \\ 0 & 1.0008 \\ 0.0180 & 0 \end{bmatrix}$$

In equation (15), the matrix $\mathbf{G}_{1,4}$ contributes to the term $(\mathbf{G}_{1,4} \otimes \mathbf{I}_N)(\Delta \mathbf{z} \otimes \Delta \mathbf{x}) + (\mathbf{I}_N \otimes \mathbf{G}_{1,4})(\Delta \mathbf{x} \otimes \Delta \mathbf{z})$. Evaluating this expression gives:

$$(\mathbf{G}_{1,4}\otimes\mathbf{I}_N)(\Delta\mathbf{z}\otimes\Delta\mathbf{x}) = \begin{bmatrix} -0.38 & 0 & 0 & 0.137 & 0 & 0 \\ 0 & -0.38 & 0 & 0 & 0.137 & 0 \\ 0 & 0 & -0.38 & 0 & 0 & 0.137 \\ 0 & 0 & 0 & 1.0 & 0 & 0 \\ 0 & 0 & 0 & 0 & 1.0 & 0 \\ 0 & 0 & 0 & 0 & 0 & 1.0 \\ 0.018 & 0 & 0 & 0 & 0 & 0 \\ 0 & 0.018 & 0 & 0 & 0 & 0 \\ 0 & 0 & 0.018 & 0 & 0 & 0 \end{bmatrix} \begin{bmatrix} \Delta x_1 \Delta z_1 \\ \Delta x_2 \Delta z_1 \\ \Delta x_3 \Delta z_1 \\ \Delta x_1 \Delta z_2 \\ \Delta x_2 \Delta z_2 \\ \Delta x_3 \Delta z_2 \end{bmatrix}$$

(82)

To transform (82) into the term $\widetilde{\mathbf{G}}_{1,4}(\Delta\mathbf{x}\otimes\Delta\mathbf{z})$, a reordering of the columns in the matrix $(\mathbf{G}_{1,4}\otimes\mathbf{I}_N)$ is needed. The reordered matrix $\widetilde{\mathbf{G}}_{1,4}$ becomes:

$$\widetilde{\mathbf{G}}_{1,4}(\Delta\mathbf{x}\otimes\Delta\mathbf{z}) \begin{bmatrix} -0.38 & 0.137 & 0 & 0 & 0 & 0 \\ 0 & 0 & -0.38 & 0.137 & 0 & 0 \\ 0 & 0 & 0 & 0 & -0.38 & 0.137 \\ 0 & 1.0 & 0 & 0 & 0 & 0 \\ 0 & 0 & 0 & 1.0 & 0 & 0 \\ 0 & 0 & 0 & 0 & 0 & 1.0 \\ 0.018 & 0 & 0 & 0 & 0 & 0 \\ 0 & 0 & 0.018 & 0 & 0 & 0 \\ 0 & 0 & 0 & 0 & 0.018 & 0 \end{bmatrix} \begin{bmatrix} \Delta x_1 \Delta z_1 \\ \Delta x_1 \Delta z_2 \\ \Delta x_2 \Delta z_1 \\ \Delta x_2 \Delta z_2 \\ \Delta x_3 \Delta z_1 \\ \Delta x_3 \Delta z_2 \end{bmatrix}$$ (83)

A similar reordering process is required when computing other coefficient matrices presented in Sections 3.2 and 3.4.

Upon examining equations (79) and (80), it becomes evident that the system exhibits discontinuities at the points $x_1 = -1$, $x_1 + z_2 = \pm\pi/2$, and $x_3 = -1$. Consequently, the system's evolution must not traverse these singularities in order for the function $\widetilde{h}$ to be well-defined and to ensure that a higher-order linearization can accurately approximate the nonlinear system. Furthermore, for the higher-order linearization to remain valid, the system trajectories must evolve within the domain: $-\pi < x_1 z_2 \leq \pi$, $-\pi < x_3 \leq \pi$, and $-\pi/2 < x_1 + z_2 \leq \pi/2$ for the higher-order linearization to be valid.

Let us now consider an initial perturbation $\Delta\mathbf{x}_0 = [-0.05 \ -0.05 \ -0.05]^T$. The complete nonlinear system response is depicted in Fig. 1, where the response is between the limits established above.

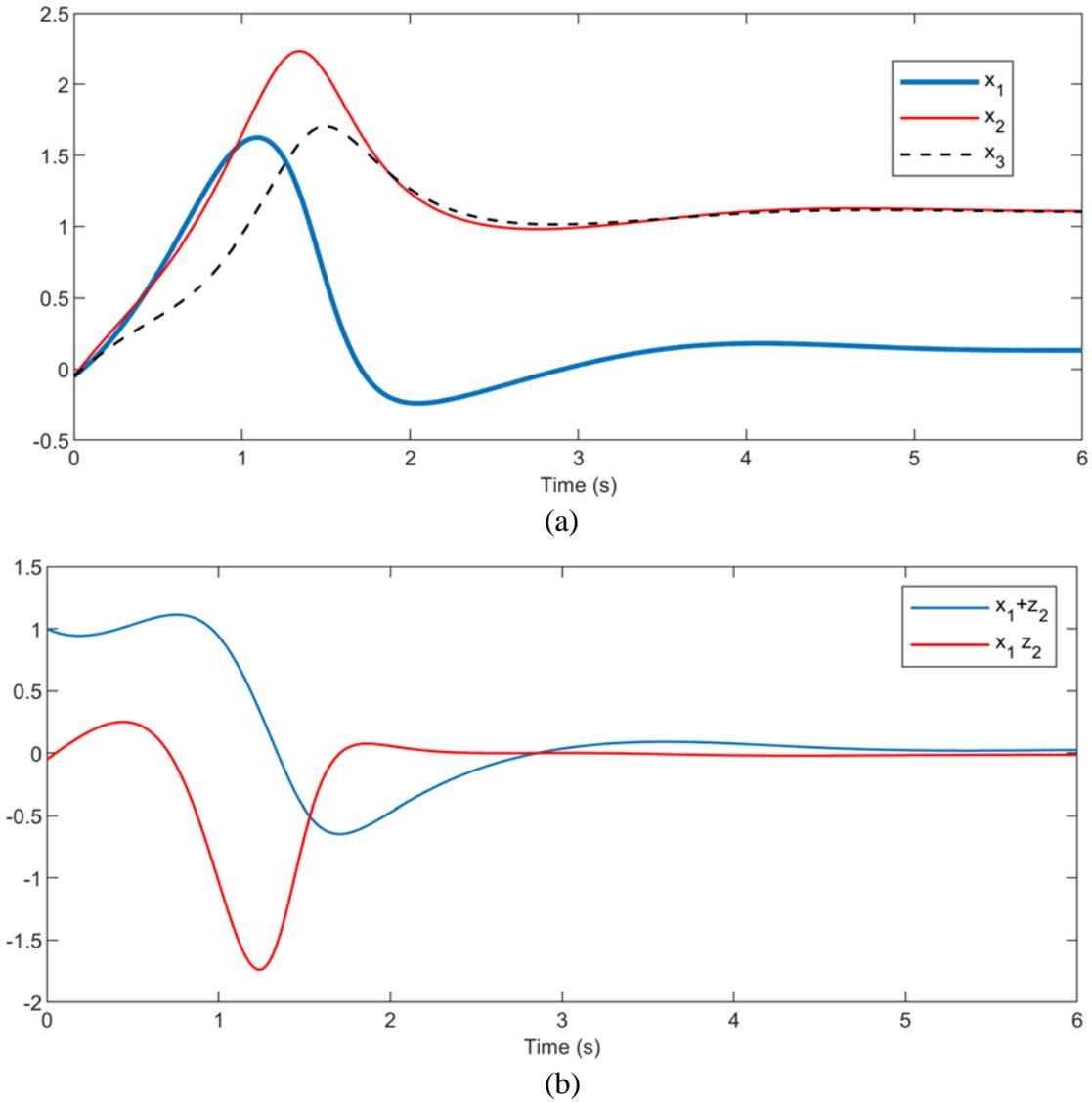

Fig. 1. Time evolution of the second test system. (a) State variables. (b) functions $x_1 + z_2$ and $x_1 z_2$.

Within this framework, the equivalent ODE system corresponding to the second- and third-order approximations was derived using equation (63). Table 1 presents the percentage difference between the matrices $\tilde{\mathbf{F}}_{11}$ and $\mathbf{A}_{nord}$, demonstrating that the proposed methodology yields a higher-order ODE system that accurately represents the original system in the vicinity of the stable equilibrium point.

The error was computed using the following expression:

$$Error = \frac{\|\tilde{\mathbf{F}}_{11} - \mathbf{A}_{nord}\|}{\|\mathbf{A}_{nord}\|} \times 100\% \qquad (84)$$

It must be emphasized that, to avoid inconsistencies arising from equivalent terms, such $\Delta x_j \Delta x_k = \Delta x_k \Delta x_j$ (a phenomenon observed in Section 4.1), the information contained in the matrices $\mathbf{A}_{nord}$ and $\tilde{\mathbf{F}}_{11}$ was preprocessed and condensed to eliminate redundancy due to such symmetric terms.

Table 1: Error of approximation of ODE matrix $\mathbf{A}_{nord}$. Second test system.

| Order of Approximation | Error (%) |
|---|---|
| Second | $1.2235 \times 10^{-15}$ |
| Third | $1.1965 \times 10^{-15}$ |

It is evident that the proposed methodology effectively captures the equivalent ODE representation of the linearized system.

*4.3 Third Test System*

As a final test case, we now consider the system defined by the following equations:

$$\begin{aligned} \dot{x}_1 &= -x_1 + \cos(x_2 - z_2) \\ \dot{x}_2 &= x_1 x_3 - x_2^2 + \sin(x_2 + z_2) \\ \dot{x}_3 &= e^{-z_1} + x_2^2 \end{aligned} \quad (85)$$

$$\begin{aligned} 0 &= x_3^2 z_1 z_2 + \cos(x_2 + z_2) \\ 0 &= -x_1 z_1^2 + z_1 x_1^2 - \sin(x_2 - z_2) \end{aligned} \quad (86)$$

In this case, the function $\tilde{h}$ is not explicitly available and is not straightforward to compute. As an initial step, we consider the first-order approximation of the system defined by equations (85)–(86), and compute the corresponding Jacobian matrix $\mathbf{A}_{1,1}$, as $\mathbf{A}_{1,1} = \mathbf{G}_{1,1} - \mathbf{G}_{1,4}(\mathbf{H}_{1,4})^{-1}\mathbf{H}_{1,1}$. The eigenvalues of $\mathbf{A}_{1,1}$ represent the linear eigenvalues of the system and are presented in Table 2.

Table 2: Linear eigenvalues of the third test system.

| Linear Eigenvalue | Frequency (Hz) | Damping Factor |
|---|---|---|
| $\lambda_1 = -1.0708$ | - | 100.0% |
| $\lambda_{2,3} = -0.133 \pm j1.7165$ | 0.2732 | 7.7266% |

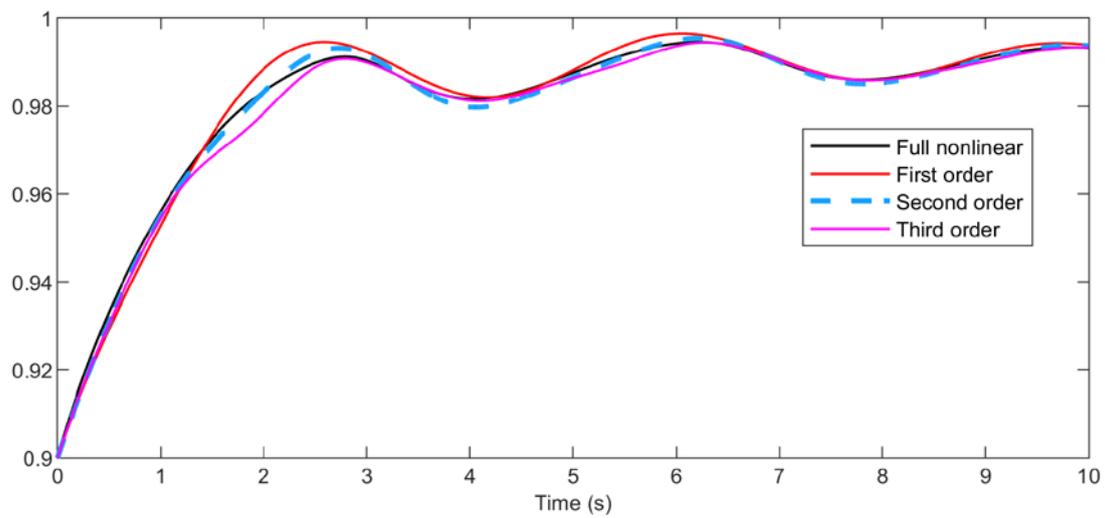

(a)

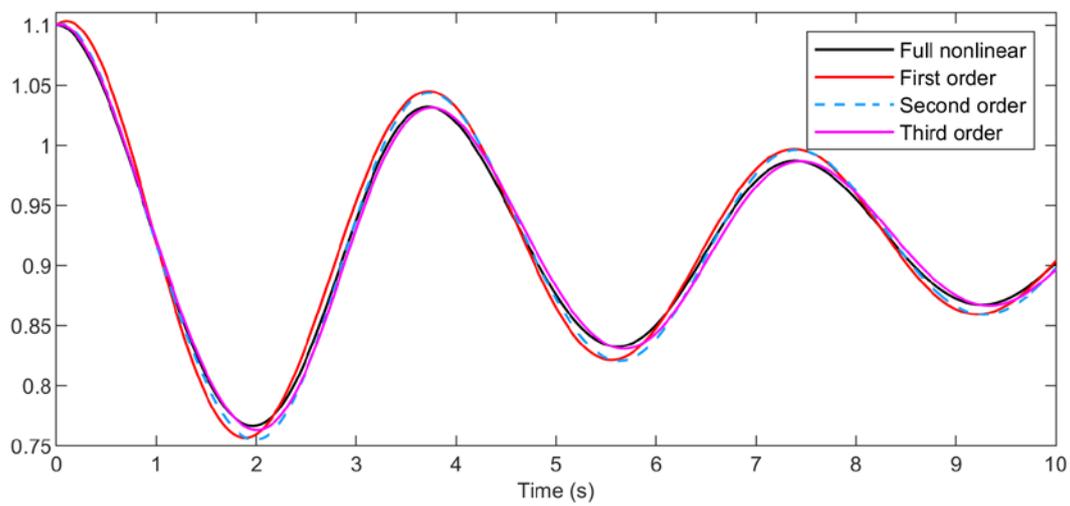

(b)

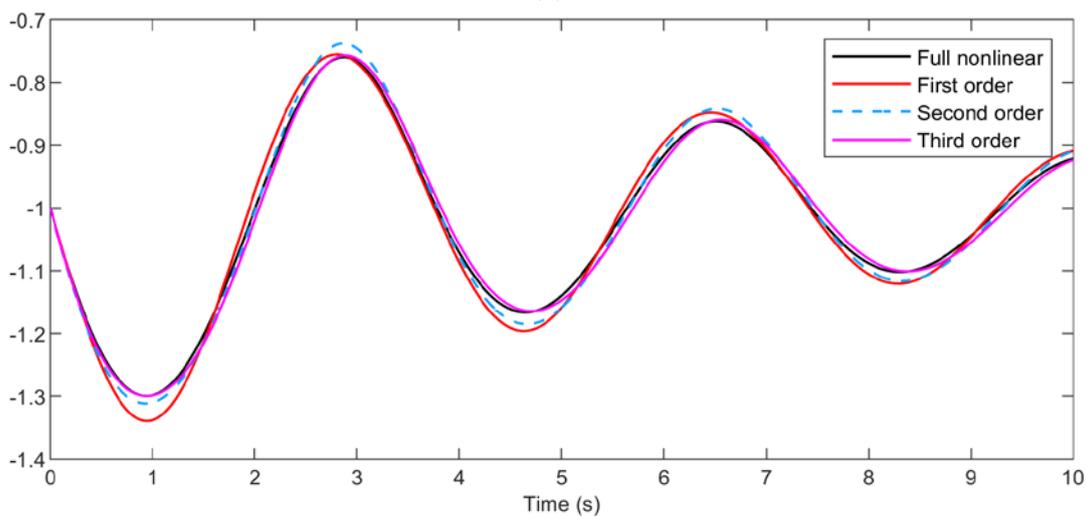

(c)

*Fig. 2. Comparison of the full nonlinear response of the third test system with the linear, quadratic, and cubic equivalent ODE systems. (a) State variable 1. (b) State variable 2. (c) State variable 3.*

Now, the linear, quadratic, and cubic extended Carleman models are constructed for this test system, and the corresponding higher-order ODE models of the form (63) are derived. To demonstrate the accuracy of these models, Fig. 2 shows the full nonlinear behavior of the system, reproduced and compared with the approximations of the linear, quadratic, and cubic models.

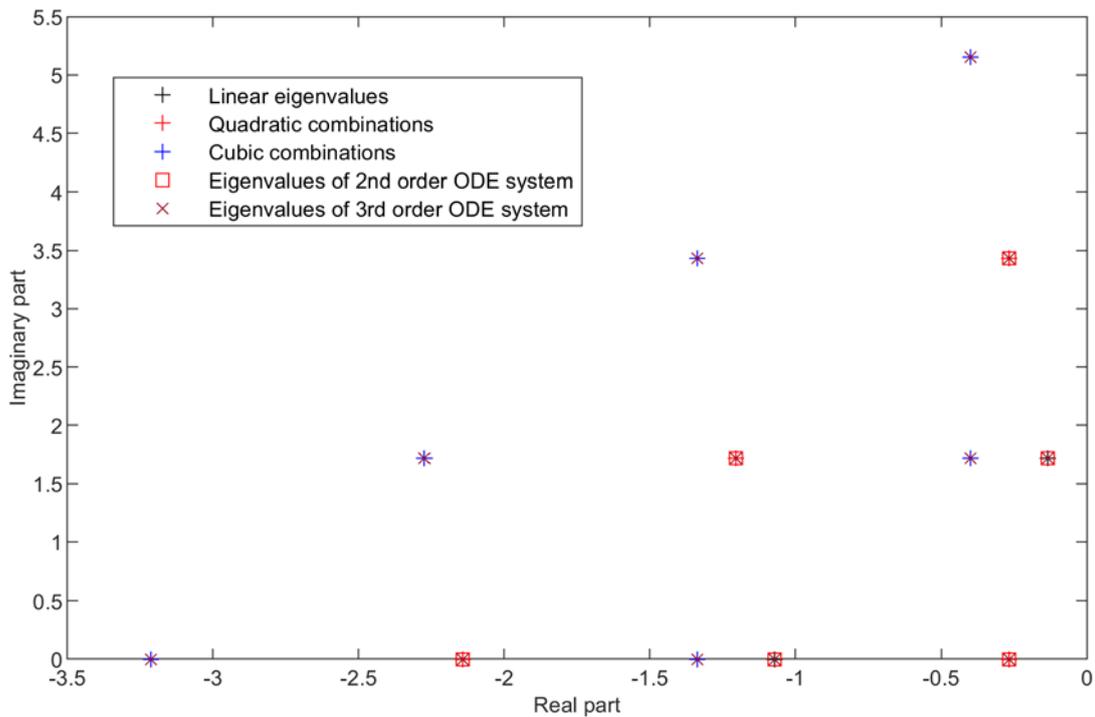

*Fig. 3. Comparison of the linear eigenvalues of the third test system and their combinations $\lambda_i + \lambda_j$ and $\lambda_i + \lambda_j + \lambda_k$ against the eigenvalues of the reduced ODE model for the second and third-order approximation.*

As shown in the results, the approximate responses become increasingly similar to the full nonlinear response as the model order increases, especially for the second and third state variables. To verify the correctness of the equivalent ODE representation, one can analyze the matrix $\tilde{\mathbf{F}}_{11}$ obtained from the quadratic and cubic approximations. This matrix should contain not only the original linear eigenvalues $\lambda_i$ but also their second- and third-order combinations, namely $\lambda_i + \lambda_j$ and $\lambda_i + \lambda_j + \lambda_k$, respectively. In Figure 3, these eigenvalues and their combinations are plotted in the complex plane and compared with the eigenvalues of $\widetilde{\mathbf{F}}_{11}$, as derived from the proposed methodology, for both second- and third-order approximations.

It can be noted that the proposed methodology successfully derives an equivalent higher-order ODE system that approximates the nonlinear response and shares the same eigenvalue combination properties as the conventional Carleman extended models.

## 5. Conclusions

This paper presented a novel methodology that extends the Carleman linearization technique to nonlinear differential-algebraic equation (DAE) systems. The conditions under which the method is applicable were explicitly stated. Subsequently, the conventional Carleman linearization rules were applied to a nonlinear DAE system.

To construct a square Carleman-extended model, additional functions derived from the algebraic equations were introduced. These functions enabled the projection of the DAE system into an equivalent higher-order ordinary differential equation (ODE) representation.

The methodology was validated through three test cases. In the first case, implementation aspects of the method were demonstrated. In the second case, the approach was shown to be effective for locally linearizable systems. It was also demonstrated that the resulting ODE representation closely matches the high-order Carleman ODE system, provided a projection from algebraic variables to state variables is available.

Finally, in the third case, where such a projection was not available, the proposed methodology still yielded an equivalent higher-order ODE model that accurately reproduced the system's nonlinear response. The resulting model preserved the key properties of Carleman linearization for ODE systems.

Future research will explore extensions of this framework, including the application of nonlinear analysis techniques such as the method of normal forms and perturbed Koopman mode analysis.